\documentclass[11pt]{article}
\usepackage{epsfig,wrapfig,epic,rotating}
\usepackage{amsmath,amsfonts,bbm,mathrsfs,stmaryrd}
\usepackage{amscd}
\usepackage{color}
\usepackage{wrapfig}
\usepackage{hyperref}

\title{Multiple Zeros of Nonlinear Systems}
\author{
Barry H. Dayton\thanks{Department of Mathematics,
Northeastern Illinois University, Chicago, IL 60625
(email: {\tt bhdayton@neiu.edu}).}
\and
 Tien-Yien Li\thanks{Department of Mathematics, Michigan State University, East Lansing, MI 48824
({\tt li@math.\allowbreak msu.edu}). Research supported in part by NSF under
Grant DMS-0811172.}
\and Zhonggang Zeng\thanks{Department of Mathematics, Northeastern Illinois
University, Chicago, IL 60625
({\tt zzeng@\allowbreak neiu.edu}). Research supported in part by NSF under
Grant DMS-0715127.}
}

\DeclareMathSymbol{\I}{\mathbin}{AMSb}{'111}
\DeclareMathSymbol{\N}{\mathbin}{AMSb}{'116}
\DeclareMathSymbol{\subsetne}{\mathbin}{AMSb}{'050}
\DeclareMathAlphabet{\mathpzc}{OT1}{pzc}{m}{it}

\newcommand{\bdb}{\mathbf{b}}
\newcommand{\bdd}{\mathbf{d}}
\newcommand{\bde}{\mathbf{e}}
\newcommand{\bdf}{\mathbf{f}}
\newcommand{\bdg}{\mathbf{g}}
\newcommand{\bdh}{\mathbf{h}}
\newcommand{\bdp}{\mathbf{p}}
\newcommand{\bdx}{\mathbf{x}}
\newcommand{\bdy}{\mathbf{y}}
\newcommand{\bdz}{\mathbf{z}}
\newcommand{\bdo}{\mathbf{0}}
\newcommand{\bdr}{\mathbf{r}}
\newcommand{\h}{{{\mbox{\tiny $\mathsf{H}$}}}}
\newcommand{\mns}{\mbox{\footnotesize -}}
\newcommand{\pls}{\mbox{\raisebox{-.4ex}{\tiny $^+$\hspace{-0.3mm}}}}
\newcommand{\npls}{\!+\!}
\newcommand{\nmns}{\!-\!}
\newcommand{\rd}{{\partial}}
\newcommand{\C}{\mathbbm{C}}

\newcommand{\bdj}{\mathbf{j}}
\newcommand{\hbdx}{{\hat{\mathbf{x}}}}
\newcommand{\hbdy}{{\hat{\mathbf{y}}}}
\newcommand{\hbdz}{{\hat{\mathbf{z}}}}
\newcommand{\sg}{\sigma}
\newcommand{\cD}{{\cal D}}

\newcommand{\dm}{\mathpzc{dim}}
\newcommand{\jet}{\mathpzc{jet}}

\newcommand{\prf}{\noindent {\bf Proof. \ }}
\newcommand{\qed}{${~} $ \hfill \raisebox{-0.3ex}{\LARGE $\Box$}}
\newcommand{\al}{\alpha}
\newcommand{\bt}{\beta}

\newcommand{\dl}{\delta}
\newcommand{\eps}{\varepsilon}
\newcommand{\spn}{\mathpzc{span}}
\newcommand{\bdk}{\mathbf{k}}
\newcommand{\bdc}{\mathbf{c}}
\newcommand{\bdi}{\mathbf{i}}
\newcommand{\cK}{{\cal K}}
\newcommand{\nullity}[1]{\mathpzc{nullity}\left(\,#1\,\right)}
\newcommand{\tms}{\mbox{\raisebox{-.4ex}{\tiny $^\times$\hspace{-0.4mm}}}}
\newcommand{\rank}[1]{\mathpzc{rank}\left(\,#1\,\right)}
\newcommand{\ranka}[2]{\mathpzc{rank}_{#1}\left(\,#2\,\right)}
\newcommand{\bdu}{\mathbf{u}}
\newcommand{\bdv}{\mathbf{v}}
\newcommand{\bdw}{\mathbf{w}}
\newcommand{\vsg}{\varsigma}
\newcommand{\pd}[2]{\frac{\partial #1}{\partial #2}}
\newcommand{\Dl}{\Delta}
\newcommand{\fM}{\mathfrak{M}}

\newtheorem{defn}{Definition}
\newtheorem{lem}{Lemma}
\newtheorem{thm}{Theorem}
\newtheorem{cor}{Corollary}
\newtheorem{example}{Example}
\newtheorem{conj}{Conjecture}

\newcommand{\aN}[2]{{\cal K}_{#1}\left(\,#2\,\right)}

\newcommand{\blb}{\big[\,}
\newcommand{\brb}{\, \big]}
\newcommand{\fm}{\mathfrak{m}}

\newcommand{\cU}{{\cal U}}

\newcommand{\hbdu}{{\hat{\mathbf{u}}}}

\newcommand{\calR}{{\cal R}}
\newcommand{\rmh}{\mathrm{h}}
\newcommand{\dd}[1]{\nabla_{#1}}

\begin{document}


%
%


\maketitle

\begin{abstract}
As an attempt to bridge between numerical analysis and algebraic geometry, 
this paper formulates the multiplicity for the general nonlinear system
at an isolated zero, presents an algorithm for computing the
multiplicity structure, proposes a depth-deflation method for
accurate computation of multiple zeros, and introduces the basic
algebraic theory of the multiplicity.
~Furthermore, this paper elaborates and proves some fundamental properties
of the multiplicity, including local finiteness, consistency, 
perturbation invarance, and depth-deflatability.
~As a justification of this formulation, the multiplicity is proved to be
consistent with the multiplicity defined in algebraic geometry for
the special case of polynomial systems.
~The proposed algorithms can accurately compute the multiplicity and
the multiple zeros using floating point arithmetic
even if the nonlinear system is perturbed.
\end{abstract}



\parskip4mm
\section{Introduction}
\vspace{-4mm}

Solving a system of nonlinear equations in the form ~$\bdf(\bdx) \,=\,\bdo$,
~or
\begin{equation} \label{nlsys}
f_1(x_1,\ldots,x_s) ~~=~~ f_2(x_1,\ldots,x_s) ~~=~~ 
\cdots ~~=~~ f_t(x_1,\ldots,x_s) ~~=~~  0 
\end{equation}
with ~$\bdf = [f_1,\ldots,f_t]^\h$ ~and ~$\bdx = (x_1,\ldots,x_s)$,
~is one of the most fundamental
problems in scientific computing, and one of the
main topics in most numerical analysis textbooks.
~In the literature outside of algebraic geometry, however, an important 
question as well as its answer seem to be absent over the years:
~What is the multiplicity of an isolated zero to the system and how to
identify it accurately.

For a single equation ~$f(x) =0$,
~it is well known that the multiplicity of a zero ~$x_*$ ~is ~$m$ ~if
\begin{equation} \label{unimult}
f(x_*) ~=~ f'(x_*) ~=~ \cdots ~=~ f^{(m\mns 1)}(x_*) ~~=~~ 0
~~~~\mbox{and}~~~ f^{(m)}(x_*) ~\ne~ 0.
\end{equation}
The multiplicity of a polynomial system at a zero has gone through rigorous
formulations since Newton's era \cite[pp. 127-129]{fulton} as one of the
oldest subjects of algebraic geometry.
~Nonetheless,
the standard multiplicity formulation and identification via
Gr\"obner bases for polynomial systems are somewhat limited to symbolic
computation, and largely unknown to numerical analysts.

As an attempt to bridge between algebraic geometry and numerical
analysis, we propose a rigorous formulation for the multiplicity structure
of a general nonlinear system at a zero.
~This multiplicity structure includes, rather than just a single integer for
the multiplicity, several structural invariances that are essential
in providing characteristics of the system and accurate computation of the 
zero.
~For instance, at the zero ~$\bdx_*=(0,0)$ ~of the nonlinear system
\begin{equation} \label{sys1}
\sin x_1 \cos x_1 - x_1 ~~=~~ \sin x_2 \sin^2 x_1 + x_2^4 ~~=~~ 0
\end{equation}
we shall have:

\vspace{-4mm}
\begin{itemize} \parskip-1mm
\item The {\em multiplicity} ~$m = 12$.
\item Under a small perturbation to the system (\ref{sys1}),
there is a cluster of exactly 12 zeros (counting multiplicities)
in a neighborhood of ~$\bdx_*=(0,0)$.
\item The {\em Hilbert function}
~$\{1,2,3,2,2,1,1,0,0,\cdots\}$ ~forms a partition of the multiplicity ~$12$.
\item
There exist 12 linearly independent differential operators
~$\rd_{00}, ~\rd_{10}, ~\ldots, ~\rd_{05}-\rd_{22},
~\rd_{06}-\rd_{23}$, grouped by the differential orders and counted by the
Hilbert function as shown in Figure~\ref{dualplot} below.
~They induce 12 differential functionals that span the
{\em dual space} associated with
system (\ref{sys1}).
~These functionals satisfy a {\em closedness condition} and
vanish on the two functions in (\ref{sys1}) at the zero $(0,0)$.
~Here, the differential operator
\begin{equation} \label{difnotation}
\rd_{j_1\cdots j_s} ~~\equiv~~
\rd_{x_1^{j_1}\cdots x_s^{j_s}} ~~\equiv~~
\mbox{\footnotesize $\displaystyle
\frac{1}{j_1!\cdots j_s!}\;\;
\frac{\rd^{j_1+\cdots+j_s}}{\rd x_1^{j_1} \cdots \rd x_s^{j_s}}$}
\end{equation}
of {\em order} ~$j_1+\cdots+j_s$ ~naturally induces a linear 
functional
\begin{equation} \label{mfnl}
\rd_{j_1\cdots j_s}[\bdx_*] ~:~ f ~\longrightarrow
(\rd_{j_1\cdots j_s} f)(\bdx_*)
\end{equation}
on functions ~$f$ ~whose indicated partial derivative exists
at the zero ~$\bdx_*$.
%
%
\item The {\em breadth}, or the nullity of the Jacobian at ~$\bdx_*$, ~is ~2.
\item The {\em depth}, which is the highest differential order
of the functionals at ~$\bdx_*$, ~is ~6.
\end{itemize}
\vspace{-4mm}

\begin{figure}[hb]
\begin{center}
\epsfig{figure=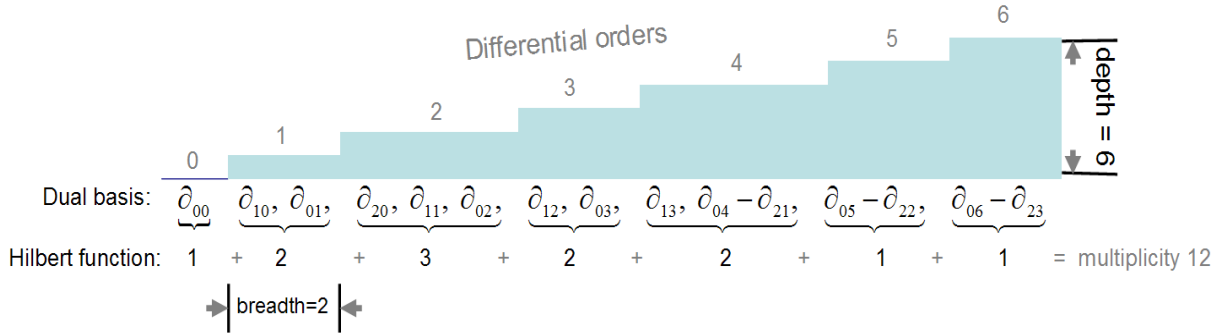,height=1.8in,width=6.3in}
\end{center} \vspace{-6mm}
\caption{\footnotesize
Illustration of the multiplicity structure
including dual basis, Hilbert function,
breadth and depth of the system (\ref{sys1}) at the zero ~$(0,0)$}
\label{dualplot}
\end{figure}

Such a multiplicity structure at an isolated zero of a general nonlinear
system will be introduced in \S\ref{s:mult}.
%
%
~We prove the so-defined {\em multiplicity} agrees with the
{\em intersection multiplicity} of polynomial systems in algebraic
geometry.
~It is finite if and only if the zero is isolated, and more importantly,
this finiteness ensures termination of the multiplicity identification
algorithm {\sc NonlinearSystemMultiplicity} given in
\S\ref{s:cms}, and it also provides
a mechanism for determining whether a zero is isolated \cite{bps09}.
~Furthermore, the multiplicity structure of the
given nonlinear system can be computed by constructing the Macaulay matrices
\cite{mac16}
together with the numerical rank revealing \cite{li-zeng-03}.
~As a result, we developed numerical algorithms that accurately calculate
the multiplicity structure even if the system data are inexact at a zero
that is given approximately (c.f. \S\ref{s:cms} and \S\ref{s:b1}).

It is well documented that multiple zeros are difficult to compute
accurately even for a single equation.
~There is a perceived barrier of ``attainable accuracy'': The number
of correct digits attainable for a multiple zero is bounded by the
number of digits in the hardware precision divided by the multiplicity.
~For instance, only three correct digits can be expected in computing a
five-fold zero using the double precision (16 digits) floating point
arithmetic.
~Such a barrier has been overcome for univariate polynomial equations
\cite{zeng-mr-05}.
~Based on the multiplicity theory established in this article,
we shall derive a {\em depth-deflation} algorithm in \S\ref{s:cmz}
for computing multiple zeros of general nonlinear systems,
which can accurately compute the multiple zeros without extending
the arithmetic precision even when the nonlinear system is perturbed.
~The {\em depth} defined in the multiplicity structure actually
bounds the number of deflation steps.
~A related {\em multiplicity deflation} method is used in
\cite{lvz06}, in which the main goal is to speed up Newton's iteration.

As mentioned above, the study of the multiplicity for a
polynomial system
at an isolated zero can be traced back to Newton's
time \cite[pp. 127-129]{fulton}.
~Besides polynomial systems, multiple zeros of a
nonlinear system occur frequently in scientific computing.
~For instance, when a system depends on certain
parameters, a multiple zero emerges when the parameters reach
a bifurcation point \cite[\S1.1]{ChoHal}.
~Accurate computation of the multiple zero and reliable identification of the
multiplicity structure may have a profound ramification in scientific
computing.
~This paper furnishes the theoretical details of
the preliminary results on polynomial systems announced
in an abstract \cite{dz}, and in addition,
the scope of this work has been
substantially expanded to general nonlinear systems.

\vspace{-4mm}
\section{Formulation and computation of the multiplicity structure}
\label{s:mult}
\vspace{-4mm}

\subsection{The notion and fundamental theorems of the multiplicity}
\label{s:not}
\vspace{-4mm}

The general nonlinear system (\ref{nlsys}) is represented by either
the mapping ~$\bdf \,:\, \C^s \longrightarrow \C^t$ ~or the set
~$F = \{ f_1, \ldots, f_t \}$ ~of functions in the variables 
~$x_1,\ldots,x_s$.
~We assume functions ~$f ~:~ \C^{s} \longrightarrow \C$ ~in this paper have
all the relevant partial
derivatives arising in the elaboration.
~The multiplicity which we shall formulate in this section
will extend both the multiplicity (\ref{unimult})
of a single equation and the Macaulay-Gr\"{o}bner duality formulation
of multiplicity for polynomial systems.

Denote ~$\N = \{0, \pm 1, \pm 2, \ldots \}$.
~For an integer array ~$\bdj = (j_1,\ldots,j_s) \in \N^s$, ~write
~$\bdj \ge 0$ ~if ~$j_i \ge 0$ ~for all ~$i\in \{1,\ldots,s\}$.
~For every ~$\bdj = (j_1,\cdots,j_s) \in \N^s$ ~with ~$\bdj \ge 0$,
~denote ~$\bdx^\bdj = x_1^{j_1}\cdots x_s^{j_s}$ ~and
~$(\bdx-\bdy)^\bdj = (x_1-y_1)^{j_1}\cdots (x_s-y_s)^{j_s}$,
~and differential functional monomial
~$\rd_\bdj[\hbdx]$
~at ~$\hbdx\in \C^s$ ~as in (\ref{mfnl}),
with order ~$|\bdj| = j_1+\cdots+j_s$.
~For simplicity, we adopt the convention
\begin{equation} \label{ng}
 \rd_\bdj [\hbdx] (f) ~~\equiv~~ 0 ~~~~\mbox{for all ~$f$ ~whenever}~~~
\bdj \not\ge 0
\end{equation}
throughout this paper.
~A linear combination
~$c = c_{\bdj_1} \rd_{\bdj_1}[\hbdx] + \cdots +
c_{\bdj_k} \rd_{\bdj_k}[\hbdx]$ ~is called a {\em differential}
{\em functional},
which will produce a set of numbers
~$c(F) =  \{c(f_1),\ldots,c(f_t) \}$
~when applied to the system ~$F = \{ f_1,\ldots,f_t \}$.
~For differential functionals, the linear anti-differentiation
transformation
~$\phi_i$ ~is defined by
~$\phi_i\big(\sum_\bdj c_\bdj \rd_\bdj[\hbdx]\big) ~=~
\sum_\bdj c_\bdj \phi_i\big(\rd_\bdj[\hbdx]\big)$
~with
\begin{equation} \label{antid}
\phi_i\big(\rd_{j_1\ldots j_s}[\hbdx]\big) ~~=~~
\rd_{j'_1\ldots j'_s}[\hbdx] ~~~~~\mbox{where}~~~~~
j'_\sg ~~=~~ \left\{ \begin{array}{cc} j_i & \mbox{if ~$\sg \ne i$} \\
j_i \nmns 1 & \mbox{if ~$\sg = i$}
\end{array} \right.
\end{equation}
for ~$i=1,\ldots,s$.
~From (\ref{ng}), we have ~$\phi_i(\rd_\bdj[\hbdx]) = 0$
~if ~$j_i = 0$.
~With these differential functionals and the linear transformations, we
now
formulate the multiplicity at a zero ~$\hbdx$ ~of the nonlinear system
(\ref{nlsys}) as follows.

\begin{defn} \label{d:mult}
~Let ~$F = \{f_1,\ldots, f_t\}$ ~be a system of functions
having derivatives of order ~$\gamma \ge 1$ ~at a zero
~$\hbdx \in \C^s$.
~Let ~$\cD_\hbdx^0 (F) \,=\,\spn\{\rd_{0\ldots 0} \}$ ~and
\begin{equation} \label{Dsg}
\mbox{$\displaystyle \cD_{\hbdx}^\al (F) \;=\;
\Big\{\, c= \sum_{\bdj \in \N^s,\,c_\bdj \in \C,\,|\bdj|\le \al}
c_\bdj \rd_\bdj[\hbdx] \;\Big|\;
c(F) = \{0\}, ~\phi_i(c) \in \cD_{\hbdx}^{\al \mns 1}(F), ~~\forall
~i=1,\dots,s \Big\}$}
\end{equation}
for ~$\al = 1, \ldots, \gamma$.
~We call such sets dual subspaces.
~If ~$\cD_\hbdx^\gamma(F) = \cD_\hbdx^{\gamma\mns 1}(F)$,
~then the vector space
\begin{equation}
\cD_\hbdx(F) ~~=~~ \cD_\hbdx^0 (F) \,\cup\, \cD_\hbdx^1 (F) \,\cup\, 
\cdots \,\cup\, \cD_\hbdx^{\gamma-1} (F) ~~=~~ \cD_\hbdx^\gamma(F)
\end{equation}
is called the \,{\em \bf dual space} of the system ~$F$ ~at ~$\hbdx$.
~The dimension of ~$\cD_\hbdx(F)$, i.e. ~$\dm \big(\cD_\hbdx(F)\big)$,
~is called the \,{\em \bf multiplicity} of ~$F$ ~at ~$\hbdx$.
\end{defn}

Notice that dual subspaces ~$\cD_\hbdx^\al(F)$'s strictly
enlarge as the differential order ~$\al$ ~increases until
reaching certain ~$\al = \dl$ ~at which ~$\cD_\hbdx^\dl(F) =  
\cD_\hbdx^{\dl\pls 1}(F)$, 
~and thus all functionals in ~$\cD_\hbdx^{\dl\pls 1}(F)$ ~are of
~differential orders up to ~$\dl$. 
~As a result, there are no functionals in the subsequent dual subspaces
with differential orders ~$\dl+2,\,\dl+3,\ldots$ ~since 
~$\phi_i\big(\cD_\hbdx^\al(F)\big) \subset \cD_\hbdx^{\al\pls 1}(F)$
~for ~$i=1,\ldots,s$. 
~Thus
\[ \cD_\hbdx^0(F) ~\subsetne~ \cD_\hbdx^1(F) ~\subsetne~ \cdots ~\subsetne
\cD_\hbdx^\dl(F) ~=~ \cD_\hbdx^{\dl\pls 1}(F) ~=~ \cdots ~=~
\cD_\hbdx^\gamma(F) ~=~ \cD_\hbdx(F).
\]
The integer ~$\dl$, ~called the {\em depth} which will be defined later,
is the highest order of differential functionals in the dual
space.

We may also denote the dual space as ~$\cD_\hbdx(\bdf)$ ~when the
nonlinear system is represented as a mapping
~$\bdf = [f_1,\ldots,f_t]^\top$.
~It is important to note that vanishing at the system ~$c(F)=\{0\}$ ~is
insufficient for the functional ~$c$ ~to be in the dual space
~$\cD_\hbdx(F)$.
~This becomes more transparent in single equation ~$f(x)=0$
~where the multiplicity is {\em not} the number of vanishing derivatives
~$f^{(k)}(x)=0$ ~at a zero ~$x_*$.
~For instance, infinite number of functionals ~$\rd_0[0], ~\rd_2[0],
~\rd_4[0], ~\ldots$ ~vanish at the ~$(1\times 1)$-system ~$\{\sin x\}$,
~since derivatives ~$\sin^{(2k)} 0 = 0$ ~for all integers ~$k\ge 0$.
~Among these functionals, however, only ~$\rd_0[0] \in \cD_0(\{\sin x\})$
~since
\[ \phi_1(\rd_{2k}[0])(\sin x) ~~=~~ \rd_{2k-1}[0](\sin x)
~~=~~ \mbox{$\frac{(-1)^{k\mns 1}}{(2k-1)!}$} \cos 0 ~~\ne~~ 0,
\]
namely ~$\rd_{2k}[0] \not\in \cD_0(\{\sin x\})$ ~for all ~$k\ge 1$,
~therefore the multiplicity of ~$\sin x$ ~is one at ~$x=0$.
~The crucial {\em closedness condition}
\begin{equation} \label{clcond}
\phi_i(c) ~~\in~~ \cD_\hbdx (F) ~~~~\mbox{for all}~~~~ c \in
\cD_\hbdx(F) ~~\mbox{and}~~ i=1,\ldots,s
\end{equation}
in Definition \ref{d:mult} requires the dual space ~$\cD_\hbdx(F)$ ~to be
invariant under the anti-differentiation transformation ~$\phi_i$'s.
~The following lemma is a direct consequence of the closedness condition.

\begin{lem} \label{l:cls}
~A differential functional ~$c$ ~is in the dual space ~$\cD_\hbdx(F)$
~of the nonlinear system ~$F = \{f_1, \ldots, f_t\}$ ~at the zero
~$\hbdx$ ~if and only if
\begin{equation} \label{cxjfi}
c\big( (\bdx-\hbdx)^\bdj f_i (\bdx) \big) ~~=~~ 0
~~~~\mbox{for any}~~~~ i\in \{1,\ldots,t\} ~~~~\mbox{and}~~~~
\bdj \in \N^s ~~~\mbox{with}~~~ \bdj \ge 0.
\end{equation}
\end{lem}

\prf
~For any ~$\bdj = (j_1,\ldots,j_s)$, ~$\bdk = (k_1,\ldots,k_s)$, ~and
function ~$f$, ~the Leibniz rule of derivatives yields
\begin{equation} \label{lr}
\rd_\bdj[\hbdx] \big( (\bdx - \hbdx)^\bdk f(\bdx) \big) ~~=~~
\rd_{\bdj\mns \bdk}[\hbdx](f)
~~\equiv~~
\big(\phi_1^{k_1}\circ \phi_2^{k_2} \circ \cdots \circ \phi_s^{k_s}\big)
(\rd_\bdj[\hbdx]) (f).
\end{equation}
The equation (\ref{cxjfi}) holds because of
the closedness condition (\ref{clcond}) and the linearity of ~$c$.
\qed

The dual space ~$\cD_\hbdx(F)$ ~itself actually contains more
structural invariants of the multiple zero beyond the multiplicity
for the system $F$.
~Via dual subspaces ~$\cD_\hbdx^\al(F)$,
~a {\em Hilbert function} ~$\rmh:\N \rightarrow \N$
~can be defined as follows:
\begin{equation}  \label{e:hildual}
\rmh(0) ~=~   \dm \big( \cD^0_{\hbdx} (F) \big) \equiv 1, ~~~~
\rmh(\al) ~=~ \dm \big(  \cD^\al_{\hbdx} (F) \big) -
\dm \big(  \cD^{\al-1}_{\hbdx} (F) \big) \;\;\mbox{\ \ for \ }
\al \in \{\,1,2,\dots \,\}.
\end{equation}
%
%
This Hilbert function is often expressed as a infinite sequence
~$\{\rmh(0),\rmh(1),\ldots \}$,
~with which we introduce the {\em breadth}
and the {\em depth} of ~$\cD_{\hbdx} (F)$,
~denoted by ~$\bt_{\hbdx}(F)$ and $\dl_{\hbdx}(F)$  ~respectively,
as 
\[ \bt_{\hbdx}(F) \; = \; \rmh\,(1) \mbox{ \ \ and \ \ }
\dl_{\hbdx}(F) \; = \; \max \{\, \al \;|\; \rmh\,(\al) > 0\, \}. \]
In other words, the breadth is the nullity of the Jacobian
at ~$\hbdx$ ~for the system (\ref{nlsys})
and the depth is the highest differential order of functionals in
~$\cD_\hbdx (F)$.
~They are important components of the multiplicity structure that dictate
the deflation process for accurate computation of the multiple zero
(c.f. \S\ref{s:cmz}).

In contrast to system (\ref{sys1}), the system
~$\{x_1^2 \sin x_1, ~x_2^2 - x_2^2 \cos x_2\}$~ also has
a zero ~$(0,0)$ ~of ~multiplicity ~$12$ ~but having a different Hilbert 
function $\{1,2,3,3,2,1,0,\cdots \}$ and a different dual space
\begin{equation} \label{db2} 
\spn\big\{\,\mbox{\footnotesize $\overbrace{\partial_{00}}^1,\;\;
\overbrace{\partial_{10} ,\;\; \partial_{01}}^2,\;\;
\overbrace{\partial_{20} ,\;\; \partial_{11}, \;\; \partial_{02}}^3, \;\;
\overbrace{\partial_{21} , \;\; \partial_{12}, \;\; \partial_{03}}^3, \;\;
\overbrace{\partial_{13} ,\;\; \partial_{22}}^2,  \;\;
\overbrace{\partial_{23}}^1$}\,\big\}.
\end{equation}
The polynomial system ~$\{x_2^3,\; x_2-x_3^2,\; x_3 - x_1^2 \}$~ at
origin is again 12-fold with Hilbert function
$\{1,\cdots,1,0,\cdots\}$ and a dual space basis
\begin{equation} \label{db3}
\mbox{\footnotesize $
\begin{array}{ll}
& \overbrace{\rd_{000}}^1,\;\;
\overbrace{\rd_{100}}^1,\;\;
\overbrace{\partial_{200}+\partial_{001}}^1, \ \cdots, \
\overbrace{\partial_{400}+\partial_{201}+\partial_{002}+\partial_{010}}^1, \\
& \cdots, \ \overbrace{\partial_{800}+\partial_{601}+\partial_{402} + \partial_{203} +
 \partial_{410} + \partial_{004} + \partial_{211} + \partial_{012} +
  \partial_{020}}^1 \\
 & \cdots, \ \overbrace{\partial_{11,00} + \partial_{901} +\partial_{702} +
   \partial_{710} + \partial_{503} +\partial_{511} + \partial_{304}
    + \partial_{312} +\partial_{105} +\partial_{320} +
    \partial_{113} +\partial_{121}}^1.
\end{array}$}
\end{equation}
The last example is of special interest
because, as a breadth-one case, its dual space can be computed via a simple
recursive algorithm (c.f. \S\ref{s:cms}).
~The dual bases in \eqref{db2} and \eqref{db3} are calculated
by applying the algorithm {\sc NonlinearSystemMultiplicity}
provided in \S\ref{s:cms} and implemented
in {\sc ApaTools} \cite{apatools}.

We now provide justifications for our multiplicity formulation
in Definition~\ref{d:mult} from its basic properties.
~First of all, the multiplicity is a direct generalization of the
multiplicity (\ref{unimult}) of univariate functions, where
the dual space at an ~$m$-fold zero ~$x_*$ ~is 
~$\cD_{x_*} (f) \,=\, \spn \{ \rd_0[x_*], ~\rd_1[x_*], ~\ldots,
~\rd_{m\mns 1}[x_*] \}$
~with Hilbert function ~$\{1,1,\ldots,1,0,\ldots\}$ ~as well as
breadth one and  depth ~$m\nmns 1$.
~Secondly, the multiplicity is well defined for analytic systems
as a finite positive integer at any isolated zero ~$\hbdx$, ~as asserted
by the Local Finiteness Theorem below.
~Thus, the process of calculating the multiplicity of an isolated zero
will always terminate at certain ~$\gamma$ ~when
~$\cD_\hbdx^\gamma(F) = \cD_\hbdx^{\gamma\mns 1}(F)$.
~The dual subspace dimensions
~$\dm\big(\cD_\hbdx^0(F)\big) \le  \dm\big(\cD_\hbdx^1(F)\big)
\le  \dm\big(\cD_\hbdx^2(F)\big) \le \cdots$
~can be unbounded if the zero lies in a higher
dimensional set of zeros.
~For example, the dual subspaces
~$\cD_{(0,0)}^\al(\{\sin(x^2), ~x\,\cos(y)\})$
~never stop expanding since infinitely many linearly independent
functionals
~$\rd_{y}[(0,0)]$, ~$\rd_{y^2}[(0,0)]$, ~$\rd_{y^3}[(0,0)]$, ~$\ldots$
~satisfy the closedness condition and vanish at the zero ~$(0,0)$.
~Obviously, ~$(0,0)$ ~lies in the zero set
~$\{(0,y)\}$, the entire ~$y$-axis,
~of the system ~$\{\sin(x^2), ~x\,\cos y \}$.

\begin{defn} \label{d:isozero}
~A point ~$\hbdx$ ~is an {\em isolated zero} of a system
~$F = \{f_1,\ldots,f_t\}$ ~if there is a neighborhood 
~$\Delta$ ~of ~$\hbdx$ ~in ~$\C^s$ 
~such that ~$\hbdx$ ~is the only zero of ~$F$ ~in ~$\Delta$.
\end{defn}

We now establish some fundamental properties of the multiplicity
for systems of {\em analytic} functions. 
~An (multivariate) analytic function, also called holomorphic function, in 
an open set ~$\Omega$ ~is commonly defined as a function ~$f$ ~that possesses
a power series expansion converging to ~$f$ ~at every point 
~$\bdx \in \Omega$ \cite[p.\,25]{Tay00}.

\begin{thm}[Local Finiteness Theorem] \label{t:lft}
~For a system ~$F$ ~of functions that are analytic
in an open set ~$\Omega \subset \C^s$,
~a zero ~$\hbdx \in \Omega$ ~is isolated if and only if
~$\sup_{\al \ge 0 } \big\{\dm \big(\cD_\hbdx^\al(F)\big)\big\}$ ~is finite. 
\end{thm}

This theorem ensures that the multiplicity is well defined at every isolated
zero, and the multiplicity computation at an isolated zero will
terminate in finitely many steps.
~It also provides a mechanism for identifying nonisolated zeros \cite{bps09}
for polynomial systems solved by homotopy method where a multiplicity 
upper bound is available. 
~The method in \cite{KuoLi} can be used to identify nonisolated zeros
for general nonlinear systems even though it is intended for polynomial
systems.

When the nonlinear system ~$P$ ~consists of polynomials
~$p_1, \ldots, p_t$ ~in the variables ~$x_1,\ldots,x_s$,
~the multiplicity theory, i.e. the {\em intersection multiplicity}
at a zero of such a special system, has been well studied in
algebraic geometry.
~The following theorem asserts that the multiplicity
~$\dm\big(\cD_\hbdx(P)\big)$ ~formulated in Definition~\ref{d:mult}
in this special case is identical to the intersection multiplicity of
polynomial systems in algebraic geometry.

\begin{thm}[Multiplicity Consistency Theorem] \label{t:2m}
~For a system ~$P$ ~of polynomials with complex coefficients,
~the multiplicity ~$\dm\big( \cD_\hbdx(P)\big)$ ~is identical to
the intersection multiplicity of ~$P$ ~at an isolated
zero ~$\hbdx$.
\end{thm}

The following Perturbation Invariance Theorem asserts that
the multiplicity as defined equals to the number of zeros ``multiplied'' from
a multiple zero when the system is perturbed.
~As a result, Definition~\ref{d:mult} is intuitively justified.

\begin{thm}[Perturbation Invariance Theorem] \label{conj}
~Let ~$F = \{f_1,\ldots,f_s\}$ ~be a system of functions that are
analytic in a neighborhood ~$\Omega$ ~of an ~$m$-fold zero 
~$\hbdx \in \C^s$ ~and ~$F^{\mns 1}(\bdo) \cap \Omega = \{\hbdx\}$.
~Then, for any functions ~$g_1,\ldots,g_s$ ~that are analytic in ~$\Omega$
~and ~$F_\eps = \{f_1 + \eps g_1, \ldots, f_s + \eps g_s\}$,
~there exists a ~$\theta >0$ ~such that, for all ~$0<\eps < \theta$,
\[
m ~=~ \dm\big(\cD_\hbdx(F)\big) ~=~ \sum_{\tilde{\bdx}\in 
F_\eps^{-1}(\bdo) \cap \Omega} 
\dm\big(\cD_{\tilde{\bdx}}(F_\eps)\big).
\]
\end{thm}

In other words, multiplicities of zeros are invariant under small perturbation
to the system of analytic functions. 
~An ~$m$-fold zero becomes a cluster of exactly ~$m$ ~zeros counting 
multiplicities.
~The proof of Theorem~\ref{conj} follows from \cite[Lemma 6]{sv}.
~We may illustrate this theorem by a computing experiment
on the following example.

\begin{example} \em
~Consider the system
~$F = \{\sin x\, \cos y - x, ~\sin y\, \sin^2 x - y^2\}$
~having multiplicity 6 at the zero ~$(0,0)$.
~In a small neighborhood of ~$(0,0)$, ~we compute the
zeros of the perturbed system
\end{example} \vspace{-3mm}
\begin{equation} \label{ptsys}
F_\epsilon ~~=~~
\{\sin x\, \cos y - x - \epsilon, ~\sin y\, \sin^2 x - y^2 + \epsilon \}
\end{equation}
for small values of ~$\epsilon$.
~A cluster of exactly 6 zeros of  ~$F_\epsilon$ ~near ~$(0,0)$ ~are
found by Newton's iteration
using zeros of the truncated Taylor series of ~$F_\epsilon$
~as the initial iterates,
matching the multiplicity of the system ~$F$ ~at $(0,0)$.
~Table~\ref{t:ptsys} shows the zeros of ~$F_\epsilon$ ~for
~$\epsilon = 10^{\mns 8}$ ~and ~$10^{\mns 12}$.
~The cluster as shown shrinks to ~$(0,0)$ ~when the perturbation
decreases in magnitude. 

\begin{table}[ht]
\begin{center}
\begin{tabular}{|c|c|} \hline
\multicolumn{2}{|c|}{$\epsilon ~~=~~ 10^{-8}$} \\ \hline
~$\bdx_1$, ~$\bdx_2$ & \scriptsize
$(-0.0039173928\mp 0.0000003908\,i, ~
-0.0000076728 \pm 0.0000997037\,i) $ \\ \hline
~$\bdx_3$, ~$\bdx_4$ &  \scriptsize
$(~~0.0019584003 \pm 0.0033883580\,i, ~
\,\,~~~0.0000035695 \pm 0.0000935115\,i) $ \\ \hline
~$\bdx_5$, ~$\bdx_6$ &  \scriptsize
$(~~0.0019590795\mp 0.0033879671\,i, ~
\,\,~~~0.0000040733\pm 0.0001067848\,i) $ \\ \hline
\multicolumn{2}{|c|}{$\epsilon ~~=~~ 10^{-12}$} \\ \hline
~$\bdx_1$, ~$\bdx_2$ & \scriptsize
$(-0.000181717560 \mp 0.000000000182\,i,~
-0.000000016511 \pm 0.000000999864\,i)$ \\ \hline
~$\bdx_3$, ~$\bdx_4$ &  \scriptsize
$(~~0.000090858627 \pm 0.000157362584\,i, ~
\,\,~~~0.000000008136 \pm 0.000000985770\,i)$ \\ \hline
~$\bdx_5$, ~$\bdx_6$ &  \scriptsize
$(~~0.000090858942 \mp 0.000157362403\,i, ~
\,\,~~~0.000000008372 \pm 0.000001014366\,i) $ \\ \hline
\end{tabular}
\end{center} \vspace{-5mm}
\caption{Zeros of the perturbed system ~$F_\epsilon$ ~in
(\ref{ptsys}) near ~$(0,0)$ ~for ~$\epsilon=10^{\mns 8}$
~and ~$10^{\mns 12}$.} \label{t:ptsys}
\end{table}

The proofs of the above three fundamental theorems on multiplicities
will be given in \S\ref{s:fou}, in which
the algebraic foundation of the multiplicity will be established.

{\bf Remark on the history of multiplicity:}
~A discussion on the history of the multiplicity formulations for
a polynomial system at a zero is given in \cite[p.127]{fulton} from
algebraic geometry.
~As Fulton points out there have been many differing concepts
about multiplicity.
~Mathematicians who have worked on this include Newton, Leibniz,
Euler, Cayley, Schubert, Salmon, Kronecker and Hilbert.
~The dual space approach was first formulated by Macaulay \cite{mac16}
in 1916 for polynomial ideals.
%
~Samuel developed this viewpoint with his {\em Characteristic functions and
polynomials} now called {\em Hilbert functions and polynomials}.
~More than the multiplicity at a zero of a polynomial system he
defines the multiplicity of an arbitrary local ring
\cite[Ch. VIII \S 10]{zs},
which, in the case of a 0-dimensional local ring, is the
sum of the Hilbert function values as in Corollary~\ref{p:oldth1}.
~As we show in \S\ref{s:fou}, this multiplicity is also the $\C$-dimension
of the local ring which is now generally accepted as the standard definition
of multiplicity in commutative algebra for isolated zeros of
systems of equations,
see Chapter 4 of \cite{clo2} for a discussion similar to that of this paper.
~Symbolic computation of {\em Gr\"obner duality} on polynomial ideals
was initiated by Marinari, Mora and M\"oller \cite{mmm96}, as well as
Mourrain \cite{bm96}.
~Stetter and Thallinger introduced numerical computation of the dual basis
for a polynomial ideal in \cite{ste-tha,Tha96}
and in Stetter's book \cite{stetter}.
~Other computational algorithms on the multiplicity problem
have recently been proposed in \cite{bps06}, \cite{kss}, \cite{lbh},
\cite{WuZhi08}, and \cite{ZengCsdual}, etc.

\vspace{-4mm}
\subsection{The Macaulay matrices} \label{s:mm}
\vspace{-4mm}

Based on the multiplicity formulation,
computing the multiplicity structure can be converted to the 
rank/kernel problem of matrices.
~Consider the dual subspace ~$\cD_\hbdx^\al (F)$ ~as defined in (\ref{Dsg})
for the nonlinear system ~$F = \{f_1,\ldots,f_t\}$ ~in ~$s \le t$
~variables ~$\bdx = (x_1,\dots, x_s)$.
%
%
~Similar to Lemma~\ref{l:cls}, one can show that a functional
~$c \,=\, \sum_{|\bdj|\le \al} c_\bdj\, \partial_{\bdj}[\hbdx]$
~is in the dual subspace ~$\cD_\hbdx^\al (F)$ ~if and only if
\begin{equation} \label{saleq}
c\big( (\bdx - \hbdx)^\bdk f_i (\bdx) \big)
~~\equiv~~ \sum_{|\bdj|\le \al} c_\bdj \cdot \partial_\bdj[\hbdx]
\big( (\bdx - \hbdx)^\bdk f_i (\bdx) \big)
~~=~~ 0
\end{equation}
for all ~$|\bdk| \le \al-1$ ~and ~$i \in \{1,\ldots,s\}$.
~By a proper ordering of indices ~$\bdj$ ~and ~$(\bdk,i)$,
~equation (\ref{saleq}) can be written in matrix form
\begin{equation} \label{salc}
S_\al \,\bdc ~~=~~ \bdo
\end{equation}
where ~$\bdc$ ~is the vector formed by ordering ~$c_\bdj$ ~in (\ref{saleq})
for ~$\bdj \in \N^s$, ~$\bdj \ge 0$ ~and ~$|\bdj|\le \al$.
~The equation (\ref{salc}) determines the dual subspace ~$\cD_\hbdx^\al (F)$
~that is naturally isomorphic
to the kernel ~$\cK(S_\al)$ ~of the matrix
~$S_\al$, ~which we call the ~$\al$-th order Macaulay matrix.

To construct the Macaulay matrices,
we choose the {\em negative degree lexicographical ordering} \cite{sing},
denoted by ~$\prec$, ~on the index set
~$\I_\al \equiv \left\{\,\bdj \in \N^s\;\big|\;\bdj \ge 0,
~|\bdj| \le \al\,\right\}$:
\begin{eqnarray*}
\lefteqn{\bdi ~\prec \bdj ~~~~\mbox{if}~~~ |\bdi| < |\bdj|, ~~~\mbox{or}}, \\
& & (|\bdi|=|\bdj| ~~~\mbox{and}~~~ \exists ~1 \le \sg \le s:
~~i_1 = j_1, \ldots, ~i_{\sg\mns 1} = j_{\sg\mns 1}, ~i_\sg < j_\sg).
\end{eqnarray*}
The Macaulay matrix ~$S_\al$ ~is of size ~$m_\al
\times n_\al$ ~where 
\[ m_\al~~=~~\left( \begin{array}{c} \al-1+s \\ \al-1  \end{array}\right) 
~~~\mbox{and}~~~
n_\al= \left( \begin{array}{c} \al+s \\ \al  \end{array}\right).
\]
We view the rows to be indexed
by ~$(\bdx-\hbdx)^\bdk\, f_i$ ~for ~$(\bdk,i) \in \I_{\al-1} \times
\{1,\cdots,t\}$ ~with ordering
~$(\bdk,i) \prec (\bdk',i')$ ~if ~$\bdk \prec \bdk'$ ~in ~$\I_{\al-1}$
~or ~$\bdk = \bdk'$ ~but ~$i<i'$, and the
columns are indexed by the differential
functionals ~$\partial_\bdj$~ for ~$\bdj \in \I_\al$.
~The entry of ~$S_\al,$ ~at the intersection of the row and column
indexed by ~$(\bdx-\hbdx)^\bdk \, f_i$ ~and ~$\partial_{\bdj}$ ~respectively,
is the value of
~$\partial_{\bdj}[\hbdx]\left((\bdx-\hbdx)^\bdk \, f_i\right)$.
~With this arrangement, ~$S_\al$ ~is the upper-left
~$m_\al \times n_\al$~ submatrix of subsequent Macaulay matrices\
~$S_\sg$,\, for \,$\sg \geq \al$, ~as illustrated in Example~\ref{mmat}.
~The following corollary is thus straightforward.

\begin{cor}  \label{p:oldth1}
~Let ~$F = \{ f_1,\dots, f_t \}$ ~be
a system of functions in variables ~$\bdx = (x_1,\dots,x_s)$
~with a zero ~$\hbdx$.
~Then for each ~$\al > 0$,
~the dual subspace ~$\cD_\hbdx^\al (F)$ ~is isomorphic to the kernel
~$\cK(S_\al)$ ~of the Macaulay matrix ~$S_\al$.
~In particular, with ~$S_0 \equiv [f_1(\hbdx),\ldots,f_t(\hbdx)]^\top = \bdo$, 
~the Hilbert function
\begin{equation} \label{e:hilnull}
\rmh(\al) ~~=~~ \nullity{S_\al} - \nullity{S_{\al\mns 1}}
~~~~\mbox{for}~~~ \al = 1, 2, \cdots.
\end{equation}
\end{cor}

Notice that for an obvious ordering ~$\prec$~ of ~$\I_1$
~and ~$\bdf(\hbdx) = [f_1(\hbdx),\ldots,f_t(\hbdx)]^\top$, ~we can arrange
\begin{equation} \label{S1}
 S_1 ~~=~~ \left[ \bdf(\hbdx) ~\big|~ J(\hbdx) \right]  \;\
~~\equiv~~ \left[ \bdo ~\big|~ J(\hbdx) \right]  
\end{equation}
where ~$J(\hbdx)$ ~is the Jacobian of the system
~$\{f_1, \dots, f_t\}$~ at \,$\hbdx$.

\begin{example} \label{mmat} {\em
~Consider the system ~$F = \{x_1-x_2+x_1^2, \;x_1-x_2+x_2^2\}$ ~at
~$\hbdx = (0,0)$.
~Figure~\ref{mmfig} shows the expansion of the Macaulay matrices from
~$S_1$ ~to ~$S_2$, ~then ~$S_3$. 
%
~The table beneath the Macaulay matrices in Figure~\ref{mmfig}
shows the bases for the kernels
as row vectors using the same column indices.
~It is instructive to compare this pair of arrays to those in
\cite[\S \,65]{mac16} or the
reconstruction of Macaulay's arrays in \cite[Example 30.4.1]{mora2}.
~For this example, the kernels can be converted to bases of dual subspaces 
using the indices in the table:
\begin{eqnarray*}
 \cD^0_{(0,0)}(F)  & = &  \spn\{\rd_{00}\}, ~~~~ 
 \cD^1_{(0,0)}(F)  ~~=~~  \spn\{\rd_{00},~~\rd_{10}+\rd_{01}\} \\
 \cD^2_{(0,0)}(F)  & = &  \spn\{\rd_{00},~~\rd_{10}+\rd_{01},
~~-\rd_{10}+ \rd_{20}+\rd_{11}+\rd_{02}\}.
\end{eqnarray*}
Since ~$\nullity{S_3} = \nullity{S_2} = 3$,
~the Hilbert function ~$ \rmh(\N) = \{1,1,1,0,\cdots \} $.
~The multiplicity equals 3.
~The dual space ~$\cD_{(0,0)}(F) = \cD^2_{(0,0)}(F)$~ with breadth
~$\beta_{(0,0)}(F) = \rmh(1) = 1$~ and depth ~$\delta_{(0,0)}(F) =
\max\{ \al \,|\, \rmh(\al)>0 \} = 2$.
~The complete multiplicity structure is in order.} \qed
\end{example}

\begin{figure}[ht]
\begin{center}
\footnotesize
$ \begin{array}{rrrr||rrrrrrrrrr|}
& \multicolumn{3}{l}{\mbox{Macaulay} }
&\multicolumn{1}{r}{\raisebox{-1.4ex}{$\overbrace{\;\;\;\;\;\;\;\;\;}^{
\mbox{\scriptsize $|\bdj|=0$}}$}} & \multicolumn{2}{c}{\raisebox{-1.4ex}{$
\overbrace{\;\;\;\;\;\;\;\;\;\;\;\;}^{
\mbox{\scriptsize $|\bdj|=1$}}$}} & \multicolumn{3}{c}{\raisebox{-1.4ex}{$
\overbrace{\;\;\;\;\;\;\;\;\;\;\;\;\;\;\;\;\;\;\;\;}^{
\mbox{\scriptsize $|\bdj|=2$}}$}} & \multicolumn{4}{c}{\raisebox{-1.4ex}{$
\overbrace{\;\;\;\;\;\;\;\;\;\;\;\;\;\;\;\;\;\;
\;\;\;\;\;\;\;\;\;\;\;}^{\mbox{\scriptsize $|\bdj|=3$}}$}}  \\
&\multicolumn{3}{r||}{\mbox{matrices $\searrow$}} &
\multicolumn{1}{r}{\rd_{00}}  & \rd_{10}  &
\multicolumn{1}{c}{\rd_{01}}  &\rd_{20}
    & \partial_{11}  & \multicolumn{1}{c}{\partial_{02}} &
\rd_{30} & \rd_{21} & \rd_{12} & \multicolumn{1}{r}{\rd_{03}} \\ \hline \hline
&& \begin{rotate}{-90}  \hspace{-0mm}
\raisebox{9.0ex}{$\underbrace{\;\;\;\;\;\;\;}_{\hspace{-3mm}
\mbox{\scriptsize $|\bdk| = 0$}}$} \end{rotate}\;\;
&&&&&&&&&&& \multicolumn{1}{c}{ } \\
&&&  f_1 & \multicolumn{1}{r|}{0} & 1 & \multicolumn{1}{r|}{-1} & 1 & 0 & 
\multicolumn{1}{r|}{0} & 0 & 0 & 0 & 0 \\
&S_0& &  f_2 & \multicolumn{1}{r|}{0} & 1 & \multicolumn{1}{r|}{-1} & 0 & 0 & 
\multicolumn{1}{r|}{1} & 0 & 0 & 0 & 0 \\ \cline{1-5}
&S_1& &  &  & & \multicolumn{1}{r|}{ } &  &  & 
\multicolumn{1}{r|}{ } &  &  &  &  \\ \cline{1-7}
&&\begin{rotate}{-90}   \hspace{-3mm}
\raisebox{6.5ex}{$\underbrace{\,\;\;\;\;\;\;\;\;\;\;\;\;\;}_{
\mbox{\scriptsize $|\bdk| = 1$}}$} \end{rotate}\;\;
&x_1f_1 &  0 & 0 &  0 & 1 & -1 & \multicolumn{1}{r|}{0} & 1 & 0 & 0 & 0 \\
&&&x_1f_2 &  0 & 0 &  0 & 1 & -1 & \multicolumn{1}{r|}{0} & 0 & 0 & 1 & 0 \\
&&&x_2f_1 &  0 & 0 &  0 & 0 & 1  & \multicolumn{1}{r|}{-1} & 0 & 1 & 0 & 0 \\
& S_2&
&x_2f_2 &  0 & 0 &  0 & 0 & 1  & \multicolumn{1}{r|}{-1} & 0 & 0 & 0 & 1 \\
\cline{1-10}
&&&x_1^2 f_1  &  0 & 0 &  0 & 0 & 0  & 0 & 1 & -1 & 0 & 0 \\
&&\begin{rotate}{-90}   \hspace{-6.0mm}
\raisebox{4.0ex}{$\underbrace{
\;\;\;\;\;\;\;\;\;\;\;\;\;\;\;\;\;\;\;\;\;}_{
\mbox{\scriptsize $|\bdk| = 2$}}$} \end{rotate}\;\;
&x_1^2 f_2  &  0 & 0 &  0 & 0 & 0  & 0 & 1 & -1 & 0 & 0 \\
&&&x_1x_2 f_1 &  0 & 0 &  0 & 0 & 0  & 0 & 0 &  1 & -1 & 0 \\
&&&x_1x_2 f_2 &  0 & 0 &  0 & 0 & 0  & 0 & 0 &  1 & -1 & 0 \\
&&&x_2^2 f_1  &  0 & 0 &  0 & 0 & 0  & 0 & 0 &  0 & 1 & -1 \\
& S_3&
&x_2^2 f_2  &  0 & 0 &  0 & 0 & 0  & 0 & 0 &  0 & 1 & -1 \\ \cline{1-14}
&&\multicolumn{2}{c}{ } &&&&&&&&&& \multicolumn{1}{c}{ }\\
\multicolumn{4}{c||}{} & \multicolumn{10}{c}{ \mbox{bases for kernels
(transposed as row vectors)} }\\
\multicolumn{4}{r||}{ } &&&&&&&&&& \multicolumn{1}{c}{ }\\
&& & \cK(S_0) & \multicolumn{1}{r|}{1} & 0 & \multicolumn{1}{r|}{0} & 0 & 0 &
\multicolumn{1}{r|}{0} &0 &0 &0 & \multicolumn{1}{r|}{0}\\ \cline{3-5}
& \multicolumn{3}{r||}{\cK(S_1)\;\;\;}
   & 0 & 1 & \multicolumn{1}{r|}{1} & 0 & 0 &
\multicolumn{1}{r|}{0} &0 &0 &0 &
\multicolumn{1}{r|}{0}\\ \cline{2-7}
\multicolumn{4}{r||}{\cK(S_2)\;\;\;\;\;\;}
 & 0 & -1 & 0 & 1 & 1 & \multicolumn{1}{r|}{1} & 0 & 0 & 0 &
\multicolumn{1}{r|}{0}\\ \cline{1-10}
\multicolumn{4}{r||}{\cK(S_3)\;\; \;\;\;\;\;\;\;\;\;}
 &   &    &   &   &   & \multicolumn{1}{r}{ } &   &   &   &
\multicolumn{1}{r|}{ }\\ \cline{1-14}
\end{array}  $
\normalsize
\end{center} \vspace{-4mm}
\caption{\footnotesize Expansion of the Macaulay matrices for the polynomial
system in Example~\ref{mmat}} \label{mmfig}
\end{figure}

By identifying the multiplicity structure of a nonlinear system
with the kernels and nullities of Macaulay matrices, the multiplicity
computation can be reliably carried out
by matrix rank-revealing, as we shall elaborate
in \S\ref{s:cms}.

\vspace{-4mm}
\subsection{Computing the multiplicity structure} \label{s:cms}
\vspace{-4mm}

The multiplicity as well as the multiplicity structure can be computed
using symbolic, symbolic-numeric or floating point computation
based on Corollary~\ref{p:oldth1}.
~The main algorithm can be outlined in the following pseudo-code.

\noindent{\bf Algorithm:} ~{\sc NonlinearSystemMultiplicity}
{\tt
\vspace{-4mm}
\begin{itemize}
\item[] Input: ~system $F = \{f_1,\cdots,f_t\}$ and isolated 
zero $\hbdx \in \C^s$
\begin{itemize}\parskip-0.5mm
\item initialize $S_0 = O_{t\tms 1}$, ~$\cK(S_0) = \spn\{[1]\}$, ~$\rmh(0) = 1$
\item for $\al = 1, 2, \cdots$ do
\begin{itemize}
\item[$*$] expand ~$S_{\al\mns 1}$ ~to ~$S_\al$, ~and
embed ~$\cK(S_{\al\mns 1})$ ~into ~$\cK(S_\al)$
\item[$*$]
find~~$\cK(S_\al)$ ~by expanding ~$\cK(S_{\al\mns 1})$
\item[$*$] if ~$\nullity{S_\al} = \nullity{S_{\al\mns 1}}$ ~then
\begin{itemize} \item[] $\dl = \al-1$, $\rmh(\al)=0$, break the loop
\item[] otherwise,
get ~$\rmh(\al)$ by (\ref{e:hilnull})
\end{itemize}
\item[] end if
\end{itemize}
\item[] end do
\item convert ~$\cK(S_\dl)$ ~to ~$\cD_{\hbdx}(F)$
\end{itemize}
\item[] Output: multiplicity ~$m = \sum_\al \rmh(\al)$, the Hilbert function
~$\rmh$, ~$\cD_{\hbdx}(F)$ ~basis, depth ~$\dl_\hbdx (F)$,
and breadth $\bt_\hbdx (F) = \rmh(1)$
\end{itemize}
}

This algorithm turns out to be essentially equivalent to
Macaulay's procedure of 1916 for finding {\em inverse arrays} of
{\em dialytic arrays} \cite{mac16,mora2}, except that
Macaulay's algorithm requires construction of dialytic arrays
with full row rank, which is somewhat difficult and costly to implement 
with inexact systems or the approximate zeros.
~Implementation of the algorithm
{\sc NonlinearSystemMultiplicity} is straightforward for
symbolic computation when the system and zero are exact and properly
represented.
~Applying this multiplicity-finding procedure on approximate zeros and/or
inexact systems requires the notions and algorithms of numerical
rank-revealing at the step ``find ~$\cK(S_\al)$'' in
Algorithm {\sc NonlinearSystemMultiplicity}.

%
%
The {\em numerical rank} of a matrix ~$A$ ~is defined as the minimum
rank of matrices within a
threshold ~$\theta$ ~\cite[\S 2.5.5]{gvl}:
~$\ranka{\theta}{A} = \min_{\|A-B\|_2 \le \theta} \rank{B}$.
The {\em numerical kernel} ~$\aN{\theta}{A}$~ of $A$ is the (exact)
kernel ~$\cK(B)$ ~of ~$B$ ~that is nearest to ~$A$ ~with
~$\rank{B}=\ranka{\theta}{A}$.
~With this reformulation, numerical rank/kernel computation
becomes well-posed.
~We refer to \cite{li-zeng-03} for details.

Numerical rank-revealing applies the iteration \cite{li-zeng-03}
\begin{equation} \label{nulvec}
\left\{ \begin{array}{l}
 \bdu_{k\pls 1} ~~=~~ \bdu_k -
\mbox{\scriptsize $
\left[ \begin{array}{c} 2\|A\|_\infty \bdu_k \\ A \end{array} \right]^\dagger
\left[ \begin{array}{c} \|A\|_\infty (\bdu_k^\h \bdu_k -1) \\ A \bdu_k
\end{array} \right] $}\\
\vsg_{k\pls 1} \, = \, \frac{\|A\bdu_{k\pls 1}\|_2}{\|\bdu_{k\pls 1}\|_2},
\;\;\;\;\; k = 0, 1, \cdots
\end{array} \right.
\end{equation}
where $(\cdot)^\dagger$ denotes the Moore-Penrose inverse.
~From a randomly chosen ~$\bdu_0$,
~this iteration virtually guarantees convergence
to a numerical null vector ~$\bdu$, ~and ~$\{\vsg_k\}$
~will converge to the distance ~$\vsg$ ~between ~$A$ ~and the nearest
rank-deficient matrix.

With a numerical null vector ~$\bdu$,
~applying (\ref{nulvec}) on
~$\hat{A} = \mbox{\scriptsize $\begin{bmatrix} \|A\|_\infty \bdu^\h
\\ A \end{bmatrix}$}$
~yields another sequence ~$\{\hat\bdu_k\}$ ~that converges to a numerical
null vector ~$\bdv$ ~of ~$A$ ~orthogonal to ~$\bdu$, ~and
the sequence ~$\{\hat\vsg_k\}$ ~converges to the distance between ~$A$
~and the nearest matrix with nullity 2.
~This process can be continued by stacking
~$\|A\|_\infty \bdv^\h$ ~on top of
~$\hat{A}$ ~and applying (\ref{nulvec}) on the new stacked matrix.

We now describe the numerical procedure for
the step of computing ~$\cK(S_\al)$  
~in Algorithm {\sc NonlinearSystemMultiplicity}.
~The kernel ~$\aN{\theta}{S_0} = \spn\{[1]\}$.
~Assume an orthonormal basis ~$Y = \blb \bdy_1,\cdots,\bdy_\mu \brb$~
for ~$\aN{\theta}{S_{\al\mns 1}}$ ~and the QR decomposition
~$ \mbox{\scriptsize $\begin{bmatrix} T Y^\h \\ S_{\al\mns 1} \end{bmatrix}$}
= Q_{\al\mns 1} \mbox{\scriptsize $\begin{bmatrix} R_{\al\mns 1} \\ O
\end{bmatrix} $}$
~are available, where ~$Q_{\al\mns 1}$ is unitary,
~$R_{\al\mns 1}$ is square
upper-triangular and ~$T$ ~is a diagonal scaling matrix.

Embedding ~$\bdy_i$'s into ~$\C^{n_\al}$ ~by appending zeros at the
bottom to form ~$\bdz_i$~ for ~$i=1,\cdots,\mu$,
~it is clear that the columns of ~$Z = \blb \bdz_1,\cdots,\bdz_\mu \brb$
~form a subset of an orthonormal basis for $\aN{\theta}{S_\al}$.
~Also, we have matrix partitions
\[
S_\al ~~=~~ \mbox{$
\begin{bmatrix} S_{\al\mns 1} & F \\ O & G \end{bmatrix}$}, ~~~~~~
\mbox{$\begin{bmatrix} T Z^\h \\ S_\al \end{bmatrix}$} ~~=~~
\mbox{
$\begin{bmatrix} T Y^\h & O \\ S_{\al\mns 1} & F \\ O & G \end{bmatrix}$}
\mbox{\scriptsize
$\begin{bmatrix} Q_{\al\mns 1}
\begin{bmatrix} R_{\al\mns 1} & F_1 \\ O & F_2 \\ \hline\end{bmatrix} \\
\;\;\;\;\;\;\;\;\;\,
\begin{bmatrix} \;\;\;O \; & \;\;\;G \; \end{bmatrix}
\end{bmatrix}$} 
\]
where ~$ \mbox{\scriptsize
$\begin{bmatrix} F_1 \\ F_2 \end{bmatrix}$} = Q_{\al\mns 1}^\h
\mbox{\scriptsize
$\begin{bmatrix} O \\ F \end{bmatrix}$}$.
%
~Let
~$\hat{Q} \mbox{\scriptsize $\begin{bmatrix} \hat{R} \\ O \end{bmatrix}$}
 = \mbox{\scriptsize
$\begin{bmatrix} F_2 \\ G \end{bmatrix}$}$~ be a QR decomposition.
~Then
\begin{equation} \label{ZSR}
 \begin{bmatrix} T Z^\h \\ S_\al \end{bmatrix} ~~=~~ Q_\al
\mbox{
$\begin{bmatrix} R_{\al\mns 1} & F_1 \\ O & \hat{R} \\
O & O  \end{bmatrix}$}  =
Q_\al \begin{bmatrix} R_\al \\ O \end{bmatrix}
\end{equation}
with a proper accumulation of ~$Q_{\al\mns 1}$ ~and ~$\hat{Q}$ ~into
~$Q_\al$.
~This implies
\[ \cK(R_\al) ~~=~~ \cK(S_\al) \bigcap \cK(Z^\h)
~~=~~ \cK(S_\al) \bigcap \aN{\theta}{S_{\al\mns 1}}^\perp. \]
Therefore ~$\aN{\theta}{R_\al}$~ consists of numerical null vectors
of $S_\al$ that are approximately orthogonal to
those of $S_{\al\mns 1}$.
~The procedure below produces the numerical kernel
$\aN{\theta}{R_\al}$.

\begin{itemize} \parskip-0.5mm
\tt
\item[$\bullet$] ~let ~$A = R_\al$
\item[$\bullet$] ~for ~$i=1,2,\cdots$~ do
\begin{itemize}\parskip-0.5mm
\item ~apply iteration (\ref{nulvec}), stop at ~$\bdu$ ~and ~$\vsg$
\newline with proper criteria
\item ~if ~$\vsg > \theta$, ~exit, end if
\item ~get ~$\bdz_{\mu+i} = \bdu$,
reset ~$A$ ~with {\scriptsize $\begin{bmatrix} \|A\|_\infty \bdu^\h
\\ A \end{bmatrix}$}
\item ~update the QR decomposition ~$A = Q R$
\end{itemize}
\item[] ~end for
\end{itemize}

Upon exit, vectors ~$\bdz_{\mu\pls 1}$, ~$\cdots$, ~$\bdz_{\mu\pls \nu}$~
are remaining basis vectors of ~$\aN{\theta}{S_\al}$ ~aside
from previously
obtained ~$\bdz_1$, ~$\cdots$, ~$\bdz_\mu$.
~Furthermore, the QR decomposition of
~{\scriptsize $\begin{bmatrix} \hat{T} \hat{Z}^\h \\ S_\al \end{bmatrix}$}
~is a by-product from a proper accumulation of orthogonal
transformations.
~Here ~$\hat{Z} = \blb \bdz_1, \cdots, \bdz_{\mu\pls \nu} \brb$ ~with a
column permutation and ~$\hat{T}$ ~is again a scaling matrix.

Algorithm {\sc NonlinearSystemMultiplicity} is implemented as a function
module in the software package {\sc ApaTools} \cite{apatools}.
~For an isolated zero of a given system along with a rank threshold,
the software produces the multiplicity, breadth, depth, Hilbert function,
and a basis for the dual space.
~The software performs symbolic (exact) computation when the rank threshold
is set to zero, and carries out numerical computation otherwise.
~An example of computing the multiplicity structure for an inexact system
at an approximate zero will be shown as Example~\ref{dbex} in \S\ref{s:dm}.

\vspace{2mm}
{\bf Remarks on computational issues:}
~For an exact system, the accuracy of a zero
~$\hbdx$ ~can be arbitrarily high using multiprecision or
a deflation method described in \S\ref{s:cmz}.
~As a result, numerical rank-revealing with sufficient low threshold
will ensure accurate multiplicity identification.
~For inexact systems, the approximate zeros may carry substantial
errors due to the inherent sensitivity.
~In this case, setting a proper threshold ~$\theta$ ~for the numerical rank
revealing may become difficult.
~The depth-deflation method given in \S\ref{s:cmz} is effective in
calculating the zeros to the highest possible accuracy that may
allow accurate identification of the multiplicity.
~However, there will always be intractable cases.
~For those systems with obtainable multiplicity structure at an approximate
solution, the rank threshold needs to be set by users according to
the magnitude of errors on the system and solution.
~Generally, the threshold should be set higher than the size of error.
%

The size increase of Macaulay matrices may become an obstacle
when the number of variables is large, compounding with high depth
~$\dl_{\hbdx}(F)$.
~Most notably, when the breadth ~$\beta_{\hbdx}(F) =1$,
~the depth will reach the maximum: ~$\dl_{\hbdx}(F) = m - 1$.
~In this situation, high order ~$\al$'s and large sizes
of ~$S_\al$ ~are inevitable.
~A special case algorithm {\sc BreadthOneMultiplicity} in \S\ref{s:b1} is
developed to deal with this difficulty.
~A recently developed closedness subspace strategy \cite{ZengCsdual}
improves the efficiency of multiplicity computation substantially
by reducing the size of the matrices. 

\vspace{-4mm}
\subsection{Proofs of Theorem~\ref{t:lft} and Theorem~\ref{t:2m}}
\label{s:fou}
\vspace{-4mm}

Theorem~\ref{t:lft} and Theorem~\ref{t:2m} are well known for zero-dimensional
polynomial systems.  
~Since a zero-dimensional system has only finitely many zeros, each zero must 
be isolated in the sense of Definition~\ref{d:isozero} 
so the content of these theorems is simply 
the classical result that ~$\dm\big( \cD_\hbdx(F)\big)$ 
~is identical to the intersection multiplicity, c.f.
\cite{G,MN,mac16}, along with 
more recent expositions by Emsalem~\cite{Em}, Mourrain~\cite{bm96} and 
Stetter~\cite{stetter}.

However these results in the case of analytic systems and nonzero-dimensional 
polynomial systems with isolated zeros are well known mainly in the folklore 
of the theory of analytic functions of several complex variables. 
~We are not aware of an explicit reference in this generality.  
~The results do follow easily, however, from the considerations of the 
last two sections and accessible facts from the
literature (e.g. \cite{Tay00}).  
~Therefore this section is a short digression sketching our proof of 
Theorems 1 and 2 and stating a few useful corollaries of these Theorems. 

We will assume in this section that ~$\hbdx = \bdo$ ~is the origin. 
~The {\em local ring} of system ~$F=\{f_1,\dots,f_t\}$ ~of analytic functions 
at ~$\bdo$ ~is ~${\cal A} = \C\{x_1,\dots,x_s\}/ F\C\{x_1,\dots,x_s\}$ ~where
~$\C\{x_1,\dots,x_s\}$ ~is the ring of all complex analytic functions in the
variables ~$x_1,\dots,x_s$ ~which converge in some neighborhood of $\bdo$ 
(c.f. \cite{clo2,Tay00}). 
~This last ring has a unique maximal ideal ~$\fM$ ~generated by
~$\{x_1,\dots,x_s\}$, ~the image of which in ~${\cal A}$ ~is the unique 
maximal ideal ~$\fm$ of $\cal A$.

We will need some notations and lemmas.
~For an analytic or polynomial function define
\begin{equation} \label{jet}
\jet(f,k) ~~=~~ \mbox{\scriptsize $\sum_{|\bdj|\le k}$} ~c_\bdj \bdx^\bdj
\end{equation}
where ~$c_\bdj\,\bdx^\bdj$ ~is the term involving ~$\bdx^\bdj$
~in the Taylor series expansion of ~$f$ ~at ~$\bdo$.
~We say that a homogeneous polynomial ~$h$ ~of total degree ~$\al$ ~is the
{\em initial form of order ~$\al$} ~of analytic or polynomial function 
~$f$ ~if ~$h=\jet(f,\al)$.

\begin{lem} \label{lftlem} \parskip0mm
~Let ~$\calR$ ~be the ring of analytic
functions on open set ~$\cU \subseteq \C^s$ ~and assume
~$\hbdx = \bdo \in \cU$.
~Let ~$F = \{f_1,\dots, f_t\} \subset \calR$ ~be a system of analytic
functions with common zero ~$\hbdx$.
~Then the following are equivalent:
\begin{itemize}\parskip0mm
\item[\em (i)] The point ~$\hbdx=\bdo \in \cU$ ~is an isolated zero of ~$F$.
\item[\em (ii)] The local ring ~$\cal A$ ~is a finite dimensional 
~$\C$-algebra.
\item[\em (iii)] There is a positive integer ~$\dl$ ~such that for all  
~$|\bdj|>\dl$ ~the monomial ~$\bdx^\bdj$ ~is the 
initial form of order ~$|\bdj|$
~of some element in ~$F\C[x_1,\dots,x_s]$.
\end{itemize}
\end{lem}

\prf 
~To prove (i) implies (ii), use R\"ukert's Nullstellensatz \cite{Tay00}
to conclude that a power of the maximal ideal ~$\fM$ ~lies in 
~$F\C\{x_1,\dots,x_s\}$, ~i.e. ~$\fm^\al =0$ ~for large ~$\al$.  
~But in  the filtration
\begin{equation} {\cal A} = \fm^0 \supseteq \fm^1 \supseteq \fm^2 \supseteq
  \dots \label{mfilt}  \end{equation}
each quotient ~$\fm^\al/\fm^{\al+1}$ ~is a ~$\C$ ~vector space of finite
dimension.  
~In this case the filtration is finite,  hence ~$\dm ({\cal A})$ 
~is finite.
   
Assuming (ii) then \eqref{mfilt} must terminate and, by Nakayama's Lemma 
\cite{Tay00},  some ~$\fm^{\dl+1} = 0$.
~Consequently ~$\bdx^\bdj \in  F\C\{x_1,\dots,x_s\}$ ~for all 
~$|\bdj| > \dl$.  
~Then each such ~$\bdx^\bdj \in  F\C\{x_1,\dots,x_s\}$
~satisfies ~$\bdx^\bdj = g_1f_1+\dots+g_tf_t$ ~for some 
~$g_1,\dots,g_t$ ~in ~$\C\{x_1,\dots,x_s\}$. 
~A straightfoward argument shows that ~$\bdx^\bdj$ ~is the initial form of
~$\jet(g_1,\al)f_1+\jet(g_2,\al)f_2+\dots + \jet(g_t,\al)f_t \in
F\C[x_1,\dots,x_s]$ ~where ~$\al=|\bdj|$, ~proving (iii).

Finally an argument using Schwartz's Lemma \cite[Exercise 4, p.35]{Tay00} 
gives (iii) implies (i). \qed 

\begin{lem} \label{IFlem} 
~The Macaulay matrix ~$S_\al$ ~of the system ~$F$ ~is row equivalent to
a matrix with linearly independent rows
\begin{equation} \label{25}
 \mbox{\scriptsize $ \left[ \begin{array}{p{1.5in}|p{.75in}}
    \rule[-.15in]{0in}{.3in}  \quad $\mathpzc{rowspace}\    S_{\al-1}$  &
 \quad\quad  $B_\al$ \\  \hline
     \rule[-.15in]{0in}{.3in} \quad \quad\quad\quad  $\bdo$ & \quad\quad
 $C_\al$ \end{array} \right]$}.
\end{equation}
Moreover, every row of the matrix block ~$C_\al$ ~can be associated with the 
intitial form of certain element of ~$F\C[x_1,\dots,x_s]$ ~by multiplying 
the entries by their
column index and adding, and these forms give a basis of the space of all
initial forms of order ~$\al$ ~on ~$F\C[x_1,\dots,x_s]$.
\end{lem}

The proof follows from the construction of ~$S_\al$.  
~We can now prove Theorem~\ref{t:lft} and Theorem~\ref{t:2m}.

{\bf Proof of Theorem~\ref{t:lft}:}  
~By Lemma~\ref{lftlem}, ~$\hbdx$ ~is an isolated zero if and only if there 
exists ~$\dl$ ~with each monomial ~$\bdx^\bdj$ ~with ~$|\bdj|>\dl$ ~being an 
initial form of some element of ~$F\C[x_1,\dots,x_s]$.  
~Since the product of a monomial and an initial form is again an initial form,
it is necessary and sufficient that all monomials of specific degree 
~$\al=\dl+1$ ~are initial forms of ~$F\C[x_1,\dots,x_s]$.  
~By Lemma~\ref{IFlem} this will happen if and only if ~$C_\al$ ~in (\ref{25})
is of full column rank.  
~This is equivalent to ~$\nullity{S_\al}=\nullity{S_{\al-1}}$ ~which by 
Corollary~\ref{p:oldth1} is equivalent to 
~$\dm(\cD^{\al-1}_\hbdx(F)) =\dm(\cD^\al_\hbdx(F))$.  
~By the closedness condition this is equivalent to
~$\dm(\cD^{\al-1}_\hbdx(F)) =\dm(\cD^\beta_\hbdx(F))$ ~for all 
~$\beta \ge \al$ ~or ~$\sup_{\al\ge 0} \dm(D^\al_\hbdx(F)) < \infty$.

{\bf Proof of Theorem \ref{t:2m}:} 
~From \eqref{mfilt}, ~$\dm( {\cal A}) = \sum_{\al=0}^\infty
\fm^\al/\fm^{\al+1}$.  
~On the other hand, from Corollary 1 and Lemma~\ref{IFlem},  
~$\dm( \cD^\al_\hbdx(F))$ ~is the sum of the dimensions of the space
of initial forms of order ~$\al$, ~$\al=0,1,\dots$.
~From the proof of \cite[Prop. 5.5.12]{singbk}, ~it follows that
~$\fm^\al/\fm^{\al+1}$ ~is isomorphic to the space of initial forms of order
~$\al$ ~and so ~$\dm( \cD^\al_\hbdx(F)) = \dm( {\cal A})$
~where ~$\cal A$ ~is the local ring of the system ~$F$ ~at ~$\hbdx=\bdo$.  
~This latter dimension is commonly known as the intersection 
multiplicity. 
\qed

Furthermore, the proof above leads to the following Depth Theorem for
an isolated zero.

\begin{cor}[Depth Theorem]  \label{thm:depth}
~Let ~$F=\{f_1,\dots,f_t\}$ ~be a system of analytic functions in an open
set of ~$\C^s$ ~at an isolated zero ~$\hbdx = \bdo$.
~Then there is a number ~$\delta = \delta_\hbdx(F)$ ~called the
{\em depth of the isolated zero ~$\hbdx$} satisfying the following
equivalent conditions.
\vspace{-0.5cm}
\begin{itemize}\parskip0mm
\item[\em (i)] ~$\delta$ ~is the highest differential order of a
functional in ~$\cD_\hbdx(F)$.
\item[\em (ii)] ~$\delta$ ~is the smallest integer so that
the Macaulay matrix ~$S_{\delta+1}$ ~is row equivalent
to a matrix
~$\Bigl[ \begin{smallmatrix}  R & B \\ 0 & C \end{smallmatrix} \Bigr] $
~where ~$C$ ~is the ~$n \times n$ ~identity matrix, where
~$n = \binom{\delta+s}{s-1}$.
\item[\em (iii)] $\dl$ ~is the smallest integer such that ~$\bdx^\bdj$ 
~is the initial form of some element of 
~$F\C[x_1,\dots,x_s]$ ~for all ~$|\bdj|>\dl$.
\end{itemize}
\end{cor}

{\bf Remark:} 
~In commutative algebra the term {\em regularity index}, {\em nil-index}  
or just {\em index} is used instead of our {\em depth}.  
~In particular the index of the ideal of the system ~$F$ ~is 
~$\dl_\hbdx(F) + 1$. 

\begin{cor} \label{corps} As in Definition~\ref{d:mult}, let 
~$F = \{f_1,\ldots, f_t\}$ ~be a system of functions 
having derivatives of order ~$\gamma \ge 1$ ~at the zero ~$\hbdx \in \C^s$. 
~If ~$\cD_\hbdx^\gamma(F) = \cD_\hbdx^{\gamma\mns 1}(F)$, ~then the polynomial
system ~$\jet(F,\gamma)$ ~has the same multiplicity structure, and hence
the same multiplicity at ~$\hbdx$ ~as ~$F$.
\end{cor}

\prf The system ~$\jet(F,\gamma)$ ~has the same Macaulay matrices up to 
~$\gamma = \delta_\hbdx(\jet(F,\gamma))$ ~as the system ~$F$ ~and hence 
~$\cD^\al_\hbdx(F)= \cD^\al_\hbdx(\jet(F,\gamma)$ ~by 
Corollary~\ref{p:oldth1}.   
\qed

Note, in particular, that this Corollary applies to any analytic system with 
an isolated zero, so such a system is locally equivalent to a polynomial 
system.

\vspace{-4mm}
\section{Accurate computation of a multiple zero by deflating its depth}
\label{s:cmz}
\vspace{-4mm}

It is well known that multiple zeros are
highly sensitive to perturbations and are therefore
difficult to compute accurately
using floating point arithmetic.
~Even for a single univariate equation ~$f(x)=0$,
~as mentioned before, there
is a perceived barrier of ``attainable accuracy'': The number of attainable
digits at a multiple zero is bounded by the hardware precision divided
by the multiplicity.
~This accuracy barrier is largely erased recently in
\cite{zeng-mr-05} for
univariate polynomial equations.
~For general nonlinear multivariate systems, we propose a general 
{\em depth-deflation} method as well as its special case variation for 
breadth one systems in this section for accurate computation of multiple
zeros without extending hardware precision even when the given system
is perturbed.

\vspace{-4mm}
\subsection{The depth-deflation method} \label{s:dm}
\vspace{-4mm}

The hypersensitivity in calculating an approximation
~$\tilde{x}_*$ ~to an ~$m$-fold zero ~$x_*$
~can be illustrated by solving ~$f(x)=x^m = 0$.
~When the function is perturbed slightly to ~$f_\eps(x) = x^m - \eps$,
~the error becomes ~$|\tilde{x}_* - x_*| \,=\, |f-f_\eps|^{\frac{1}{m}}$. 
~The asymptotic condition number is
~$\sup_{\eps > 0} \,\frac{|\tilde{x}_*-x_*|}{|f-f_\eps|} \,=\, \infty$
~when the multiplicity ~$m>1$.
~Consequently, multiple zeros are referred to as ``singular'' or
``infinitely sensitive'' to perturbations in the literature.
~On the other hand, a simple zero is considered ``regular'' with a
finite condition number as stated in the following lemma.

\begin{lem} \label{l:cond}
~Let ~$\bdf$ ~be a system of ~$s$-variate functions that are twice 
differentiable in a neighborhood of ~$\hbdx \in \C^s$.
~If the Jacobian ~$J(\hbdx)$ ~of ~$\bdf(\bdx)$ ~at ~$\hbdx$ ~is injective
so that ~$\|J(\hbdx)^+\|_2 < \infty$,
then
\begin{equation} \label{Jcond}
\big\| \tilde{\bdx} - \hbdx \big\|_2 ~~\le~~
\big\| J(\hbdx)^+ \big\|_2 \, \big\| \bdf(\tilde{\bdx}) - \bdf(\hbdx) \big\|_2
+ O\big(\| \bdf(\tilde{\bdx}) - \bdf(\hbdx) \|_2^2 \big)
\end{equation}
for ~$\tilde{\bdx}$ ~sufficiently close to ~$\hbdx$.
\end{lem}

\prf
~The injectiveness of ~$J(\hbdx)$ ~implies ~$t \ge s$ ~and
~$\rank{J(\hbdx)} = s$.
~Without loss of generality, we assume the submatrix of ~$J(\hbdx)$
~consists of its first ~$s$ ~rows is invertible.
~By the Inverse Function Theorem, the function
~$[y_1,\ldots,y_s]^\h = [f_1(\bdx),\ldots,f_s(\bdx)]^\h$
~has a continuously differentiable inverse
~$\bdx = \bdg(y_1,\ldots,y_s)$ ~in a neighborhood of
~$[\hat{y}_1,\ldots,\hat{y}_s]^\h = [f_1(\hbdx),\ldots,f_s(\hbdx)]^\h$,
~permitting ~$\|\bdx-\hbdx\|_2 \le C \|\bdf(\bdx)-\bdf(\hbdx)\|_2$
~for ~$\bdx$ ~in a neighborhood of ~$\hbdx$.
~Since
\[
\bdf(\bdx) - \bdf(\hbdx) ~~=~~ J(\hbdx)(\bdx - \hbdx) + \bdr(\bdx)
~~~\mbox{or}~~~
\bdx - \hbdx ~~=~~ J(\hbdx)^+ \big[\bdf(\bdx) - \bdf(\hbdx)
- \bdr(\bdx)\big]
\]
where ~$\|\bdr(\bdx)\|_2 = O\big(\|\bdx-\hbdx\|_2^2\big) =
O\big(\|\bdf(\bdx)-\bdf(\hbdx)\|_2^2\big)$,
~we thus have (\ref{Jcond}). \qed

In light of Lemma~\ref{l:cond}, we may define
the {\bf condition number}
of the system ~$\bdf$ ~at a zero ~$\hbdx$:
\begin{equation} \label{zcond}
\kappa_\bdf(\hbdx) ~~=~~ \left\{ \begin{array}{ccl} \|J(\hbdx)^+\|_2
& & \mbox{if ~$J(\hbdx)$ ~is injective} \\
\infty & & \mbox{otherwise.} \end{array} \right.
\end{equation}
This condition number serves as a sensitivity measurement 
in the error estimate
\begin{equation} \label{errest}
\|\tilde{\bdx} - \hbdx \|_2 ~~\approx~~ \kappa_\bdf(\tilde{\bdx}) \cdot
\|\bdf(\tilde{\bdx})\|_2
\end{equation}
of the approximate zero ~$\tilde{\bdx}$ ~using the residual
~$\|\bdf(\tilde{\bdx})\|_2$.

Solving a nonlinear system for a multiple zero is an ill-posed problem
in the sense that its condition number is infinity 
\cite[Definition 1.1, p.\,17]{DemBook}.
~The straightforward Newton's iteration attains only a few correct
digits of the zero besides losing its quadratic convergence rate,
if it converges at all.
~Similar to other ill-posed problems, accurate computation of a multiple zero
needs a regularization procedure.
~An effective regularization approach is {\em deflation}
\cite{lvz06,lvz08,ojika}.
~For instance, Leykin, Verschelde and Zhao \cite{lvz06}
propose a deflation method and a higher-order deflation method \cite{lvz08}
which successfully restore the quadratic convergence of Newton's iteration.
~From our perspective, perhaps the most important feature of deflation
strategy should reside in transforming an ill-posed zero-finding into a
well-posed least squares problem.
~As a result, the multiple zero can be calculated to high accuracy.

We hereby propose two new versions of the deflation method, 
both are refered to as {\em depth-deflation} methods,
with one for the general cases and the other for the cases where the 
breadth of the system is one at the zero.
~We first derive our general depth-deflation method here.
~The version for breadth-one systems follows in \S\ref{s:b1}.

Let ~$\bdf ~:~ \C^s \longrightarrow \C^t$ ~represent
a nonlinear system ~$\bdf(\bdx)=\bdo$ ~where
~$\bdf(\bdx) = [f_1(\bdx),\cdots,f_t(\bdx)]^\top$,
~$\bdx = (x_1,\ldots,x_s)\in \C^s$
~with ~$t\ge s$,  ~and ~$\hbdx$ ~be an isolated zero of ~$\bdf(\bdx)$.
~Denote ~$J(\bdx)$ ~as the Jacobian of ~$\bdf(\bdx)$.
~If ~$\hbdx$ ~is a simple zero, then ~$J(\hbdx)$ ~is injective with
~pseudo-inverse ~$J(\hbdx)^+ = [J(\hbdx)^\h J(\hbdx)]^{\mns 1} J(\hbdx)^\h$,
~and the Gauss-Newton iteration
\begin{equation} \label{gnitf}
\bdx^{(n\pls 1)} ~~=~~\bdx^{(n)} - J(\bdx^{(n)})^+ \,\bdf(\bdx^{(n)})
~~~~\mbox{for}~~~ n = 0, 1, \ldots
\end{equation}
locally converges to ~$\hbdx$ ~at a quadratic rate.
~More importantly in this regular case,
solving ~$\bdf(\bdx)=\bdo$ ~for the solution ~$\hbdx$ ~is
a well-posed problem and the condition number
~$\|J(\hbdx)^+\| < \infty$.

When ~$\hbdx$ ~is a multiple zero of the system ~$\bdf$, ~however,
the Jacobian ~$J(\hbdx)$ ~is rank-deficient.
~In this {\em singular} case, the zero ~$\hbdx$ ~is
{\em underdetermined} by the system ~$\bdf(\bdx)=\bdo$ ~because it is
also a solution to ~$J(\bdx)\bdy = \bdo$ ~for some ~$\bdy \ne \bdo$.
~In order to eliminate the singularity and thus to curb the hypersensitivity,
~perhaps further constraints should be imposed.

Let ~$n_1 \,=\,\nullity{J(\hbdx)}$ ~which is strictly positive at the multiple
zero ~$\hbdx$.
~Denote ~$\bdx_1 = \bdx$ ~and ~$\hbdx_1 = \hbdx$.
~Then, for almost all choices of an ~$n_1 \times s$~ random matrix ~$R_1$,
~the matrix ~{\scriptsize $\left[ \begin{array}{c} J(\hbdx_1) \\ R_1 \end{array}
\right]$}~ is of full (column) rank.
~It is easy to see that the linear system
~{\scriptsize $ \begin{bmatrix}  J(\hbdx_1) \\ R_1  \end{bmatrix}
\bdx_2 = \begin{bmatrix} \bdo \\ \bde_1 \end{bmatrix} $}~
has a unique solution ~$\bdx_2 = \hbdx_2 \ne \bdo$.
~Here ~$\bde_1$ ~is the first canonical vector ~$[1,0,\ldots,0]^\top$
~of a proper dimension.
~As a result, ~$(\hbdx_1,\hbdx_2)$ ~is an isolated zero of a new 
~$(2t+k)\times (2s) $~ system
\begin{equation} \label{dflt1}
\bdf_1(\bdx_1,\bdx_2)  ~~\equiv~~
\mbox{\scriptsize $
\left[ \begin{array}{c}
\bdf(\bdx_1) \\ \left[ \begin{array}{c} J(\bdx_1) \\ R_1 \end{array} \right]
\bdx_2 - \left[ \begin{array}{l} \bdo \\ \bde_1 \end{array} \right]
\end{array} \right]$}.
\end{equation}

If ~$(\hbdx_1,\hbdx_2)$ ~is a simple zero of ~$\bdf_1(\bdx_1,\bdx_2)$, ~then
the singularity of ~$\bdf(\bdx)$ ~at ~$\hbdx$ ~is ``deflated'' by solving
~$\bdf_1(\bdx_1,\bdx_2)=0$ ~for ~$(\hbdx_1,\hbdx_2)$ ~as a well-posed problem
using the Gauss-Newton iteration (\ref{gnitf}) on ~$\bdf_1$.
~However, ~$(\hbdx_1,\hbdx_2)$ ~may still be a multiple zero of
~$\bdf_1(\bdx_1,\bdx_2)$ ~and, in this case, 
~we can repeat the depth-deflation method above on ~$\bdf_1$.
~Generally, assume ~$(\hbdx_1,\ldots,\hbdx_{2^\al})$ ~is 
an isolated multiple zero of ~$\bdf_\al(\bdx_0,\ldots,\bdx_{2^\al})$
~after ~$\al$ ~steps of depth-deflation
~with a Jacobian ~$J_\al(\hbdx_1,\ldots,\hbdx_{2^\al})$ ~of nullity
~$n_\al>0$.
~The next depth-deflation step expands the system to
\begin{equation} \label{dflta}
\bdf_{\al\pls 1}(\bdx_1,\ldots,\bdx_{2^{\al\pls 1}})  ~~\equiv~~
\mbox{\scriptsize $
\left[ \begin{array}{c}
\bdf_\al (\bdx_1,\ldots,\bdx_{2^\al}) \\
\left[ \begin{array}{c} J_\al(\bdx_1,\ldots,\bdx_{2^\al} ) \\ R_{\al\pls 1}
\end{array}
\right]
\left[ \begin{array}{l} \bdx_{2^\al +1} \\ ~~~\vdots \\
\bdx_{2^{\al\pls 1}} \end{array} \right]
- \left[ \begin{array}{l} \bdo \\ \bde_1 \end{array} \right]
\end{array} \right]$}
\end{equation}
where ~$R_{\al\pls 1}$ ~is a randomly selected
matrix of ~$n_{\al\pls 1}$ ~rows and the
same number of columns as ~$J_\al(\bdx_1,\ldots,\bdx_{2^\al})$.
~The depth-deflation process continues by expanding
~$\bdf(\bdx_1)$ ~to ~$\bdf_1(\bdx_1,\bdx_2)$, ~$\bdf_2(\bdx_1,\ldots,\bdx_4),
\,\ldots$ until reaching an expanded
system ~$\bdf_\sg(\bdx_1,\bdx_2,\ldots,\bdx_{2^\sg})$ ~with an isolated zero
~$(\hbdx_1,\ldots,\hbdx_{2^\sg})$ ~that is no longer singular.
~The following Depth Deflation Theorem ensures the
deflation process will terminate and the number of deflation
steps is bounded by the  depth ~$\dl_\hbdx(\bdf)$.

\begin{thm}[Depth Deflation Theorem] \label{t:conj1}
~Let ~$\hbdx$ ~be an isolated zero of a system ~$\bdf$ ~with 
depth ~$\dl_\hbdx(\bdf)$.
~Then there is an integer ~$\sg \le \dl_{\hbdx}(\bdf)$ ~such that the
depth-deflation process terminates at the expanded system
~$\bdf_\sg(\bdx_1,\ldots,\bdx_{2^\sg})$ ~with a simple zero
~$(\hbdx_1,\ldots,\hbdx_{2^\sg})$ ~where ~$\hbdx_1 = \hbdx$.
~Furthermore, the depth-deflation method generates ~$2^\sg$ 
~differential functionals in the dual space ~$\cD_{\hbdx}(\bdf)$. 
\end{thm}

We shall prove this Depth Deflation Theorem
via multiplicity analysis in \S\ref{s:dtech}.

For polynomial systems, Leykin, Verschelde and Zhao
proved that each deflation step of their method deflates intersection
multiplicity by at least one \cite[Theorem 3.1]{lvz06}.
~Theorem~\ref{t:conj1} improves the deflation bound substantially since
the depth is much smaller than the multiplicity when the breath is
larger than one.
~The computing cost increases exponentially as the depth-deflation continues
since each depth-deflation step doubles the number of variables.
~Fortunately, computing experiments suggest that, for a multiple
zero of breadth larger than one, very few depth-deflation steps are
required.
~At breadth-one zeros, we shall derive a special case deflation
method in \S\ref{s:b1}.
~The high accuracy achieved by applying the depth-deflation method
can be illustrated in the following examples.

\begin{example} \label{dbex} \em
~Consider the system
\end{example}
\vspace{-4mm}
\begin{equation} \label{sys123}
\left\{ \begin{array}{rcl}
(x-1)^3+\mbox{\scriptsize.416146836547142}\,(z-3)\sin y
+\mbox{\scriptsize .909297426825682}\,(z-3)\cos y & = & 0 \\
(y-2)^3+\mbox{\scriptsize .989992496600445}\,(x-1)\sin z
+\mbox{\scriptsize .141120008059867}\,(x-1)\cos z & = & 0 \\
(z-3)^3-\mbox{\scriptsize .540302305868140}\,(y-2)\sin x
+\mbox{\scriptsize .841470984807897}\,(y-2)\cos x & = & 0
\end{array} \right.
\end{equation}
which is a perturbation of magnitude ~$10^{\mns 15}$ ~from an exact system
~$\{u^3+w\,\sin v = v^3 + u\, \sin w  = 
w^3 + v\,\sin u = 0 \}$ ~with ~$u = x-1$, ~$v=y-2$ ~and 
~$w=z-3$. 
~This system has a zero ~$(1,2,3)$ ~of multiplicity 11, depth 4 and breadth 3.
~Using 16-digit arithmetic in Maple to simulate the hardware precision,
Newton's iteration without depth-deflation attains only 4 correct digits,
whileas a single depth-deflation step eliminates the singularity and 
obtains 15 correct digits, as shown in the following table.
~The error estimates listed in the table are calculated
using (\ref{errest}) which provides an adequate accuracy measurement
for the computed zeros.

\footnotesize
\begin{center}
\begin{tabular}{|c|c||l|l||l|}\hline
\multicolumn{2}{|c||}{ }  & \multicolumn{1}{|c|}{without deflation}
& \multicolumn{1}{|c||}{with deflation}
& \multicolumn{1}{|c|}{exact value}
\\ \hline
& $x$ & ~~~1.0003 & ~0.999999999999999 & ~~~1.0\\
zero & $y$ & ~~~1.9997 & ~1.999999999999999 & ~~~2.0 \\
& $z$ & ~~~3.0003 & ~3.000000000000000 & ~~~3.0 \\ \hline
\multicolumn{2}{|c||}{error estimate} & ~~~0.00027
& ~0.000000000000019 & \\ \hline
\end{tabular}
\end{center}
\normalsize

Since the estimated error of the approximate zero is 
~$1.94 \times 10^{\mns 14}$, ~we set the rank threshold to be slightly larger:
~$10^{-12}$.
~Algorithm {\sc NonlinearSystemMultiplicity} accurately produces
the multiplicity 11, breadth 3, depth 4, Hilbert function
~$\{1,3,3,3,1,0,\ldots,\}$ ~and (approximate) dual basis
\begin{eqnarray*}
\lefteqn{\rd_{000}, ~~\rd_{100}, ~~\rd_{010}, ~~\rd_{001}, ~~\rd_{200},
~~\rd_{020}, ~~\rd_{002},
~~\mbox{\scriptsize .707106781186544}\,\rd_{101}
+\mbox{\scriptsize .707106781186543}\,\rd_{030},
}  \\
& &
~~\mbox{\scriptsize .707106781186544}\,\rd_{011}
+\mbox{\scriptsize .707106781186545}\,\rd_{300},
~~\mbox{\scriptsize .707106781186545}\,\rd_{110}
+\mbox{\scriptsize .707106781186545}\,\rd_{003}, \\
& & \mbox{\scriptsize .500000000000008}\,\rd_{111}
+\mbox{\scriptsize .500000000000007}\,\rd_{400}
+\mbox{\scriptsize .500000000000009}\,\rd_{040}
+\mbox{\scriptsize .500000000000008}\,\rd_{004}.
\end{eqnarray*}

\begin{example} \em
~Consider the system
\end{example}
\vspace{-8mm}
\begin{eqnarray*}
e^z-\mbox{\scriptsize .944956946314738}\,\cos y
+\mbox{\scriptsize .327194696796152}\,\sin y  & = & 0\\
z^2-y^3-y^2-\mbox{\scriptsize .333333333333333}\,y
-\mbox{\scriptsize .0370370370370370} & = & 0 \\
y^2+\mbox{\scriptsize .666666666666667}\,y
+\mbox{\scriptsize .148148148148148}-x^3+x^2-
\mbox{\scriptsize .333333333333333}\,x & = & 0.
\end{eqnarray*}
This is a perturbation of magnitude ~$10^{\mns 15}$ ~from an exact system
~$\big\{e^z -\cos\big(y+\frac{1}{3}\big) = z^2-(y+\frac{1}{3}\big)^3 =
(y+\frac{1}{3}\big)^2 - (x-\frac{1}{3}\big)^3 = 0\big\}$ 
with zero ~$(1/3,-1/3,0)$ ~of multiplicity 9, depth 5, breadth 2 and Hilbert
function ~$\{1,2,2,2,1,1,0,\ldots \}$.
~Again, using 16-digits arithmetic in Maple, Newton's iteration diverges
from the initial iterate ~$(0.31,-0.31,0.01)$.
~In contrast, our depth-deflation method takes three deflation steps
to eliminate the singularity and obtains 15 correct digits of the multiple
zero:

\footnotesize 
\begin{center}
\begin{tabular}{|c|c||l|r||r|}\hline
\multicolumn{2}{|c||}{ }  & \multicolumn{1}{|c|}{without deflation}
& \multicolumn{1}{|c||}{with deflation}
& \multicolumn{1}{|c|}{exact value}
\\ \hline
& $x$ & ~~diverges      & 0.3333333333333336 &   $1/3$~~~~ \\
zero & $y$ & ~~diverges & -0.3333333333333334 & $-1/3$~~~~ \\
& $z$ & ~~diverges      & 0.0000000000000002 &   $0$~~~~ \\ \hline
\multicolumn{2}{|c||}{error estimate} & ~~~-----
& ~0.0000000000001950 & \\ \hline
\end{tabular}
\end{center}
\normalsize

\vspace{-4mm}
\subsection{Multiplicity analysis of the depth-deflation method} \label{s:dtech}
\vspace{-4mm}

We shall use some additional differential notations and operations.
~The original variables ~$\bdx = [ x_1,\cdots,x_s]^\top$ ~will be
denoted by ~$\bdx_1$ ~in accordance with the notation for
the auxiliary (vector) variables ~$\bdx_2$, ~$\bdx_3,\,\ldots$ ~etc.
~For any fixed or variable vector ~$\bdy = [y_1,\cdots,y_s]^\top$,
~the {\em directional differentiation operator} along ~$\bdy$ ~is defined as
\begin{equation} \label{nabla}
 \dd{\bdy} ~~\equiv~~ \mbox{$y_1 \pd{}{x_1} + \cdots +
y_s \pd{}{x_s}$}.
\end{equation}
When ~$\bdy$ ~is fixed in ~$\C^s$, ~$\dd{\bdy}$
induces a functional ~$\dd{\bdy}[\hbdx]\,:\, p \longrightarrow
(\dd{\bdy}p)(\hbdx)$.
~For any variable ~$\bdu = [u_1,\cdots,u_s]^\top$,  ~the
{\em gradient operator} ~$ \Dl_{\bdu} \equiv \left[\,\pd{}{u_1},\; \cdots,\;
\pd{}{u_s}\,\right]^\top$, ~whose ``dot product'' with a vector
~$\bdv = [v_1,\cdots,v_s]^\top$ ~is defined as
\begin{equation} \label{Delta}
 \bdv \cdot \Dl_\bdu ~~\equiv~~   \mbox{$v_1 \pd{}{u_1} + \cdots +
v_s \pd{}{u_s}$}. \end{equation}
In particular, ~$\dd{\bdy} \equiv \bdy \cdot \Dl_\bdx \equiv \bdy \cdot
\Dl_{\bdx_1}$~ for any ~$\bdy$ of dimension ~$s$.
~Let ~$\bdy$ ~and ~$\bdz$ ~be auxiliary variables.
~Then, for any function ~$f(\bdx)$,
\begin{equation}
\begin{array}{l}
\mbox{$(\bdy \cdot \Dl_{\bdx_1} ) (\dd{\bdz} f (\bdx_1))$}
~~=~~
\mbox{$\dd{\bdy} \dd{\bdz} f(\bdx_1)$},  \;\;\;\;\;\;
\mbox{
$\bdz \cdot \Dl_{\bdy} f(\bdx_1) ~~\equiv~~ 0$},
\\
~~~~\mbox{$
(\bdz \cdot \Dl_\bdy ) (\dd{\bdy} f (\bdx_1)) ~~=~~
(\bdz \cdot \Dl_\bdy ) (\bdy \cdot \Dl_{\bdx_1} ) f(\bdx_1)
 = \dd{\bdz} f(\bdx_1)$}.
\end{array}
\label{drep}
\end{equation}
Let ~$\bdf_0(\bdx_1) \equiv \bdf(\bdx) = [f_1(\bdx),\cdots,f_t(\bdx)]^\top$
~be a nonlinear system in variable vector ~$\bdx$ ~and ~$J_0(\bdx)$
~be its Jacobian matrix.
~Then
\[ J_0(\bdx) \,\bdz ~~=~~ \mbox{\scriptsize
$\begin{bmatrix} \Dl_{\bdx} f_1(\bdx)^\top \\ \vdots
\\ \Dl_{\bdx} f_t(\bdx)^\top \end{bmatrix}$}\, \bdz ~~=~~
\mbox{\scriptsize 
$\begin{bmatrix} \bdz \cdot \Dl_{\bdx} f_1(\bdx) \\ \vdots
\\ \bdz \cdot \Dl_{\bdx} f_t(\bdx) \end{bmatrix}$}
~~=~~ \dd{\bdz} \bdf(\bdx_1).
\]

The first depth-deflation step expands the system to 
~$\bdf_1(\bdx_1,\bdx_2)=\bdo$ ~with
\begin{equation} \label{lvzdfsys2}
\bdf_1(\bdx_1, \bdx_2)  ~~\equiv~~
\mbox{\scriptsize $
\left[ \begin{array}{c}  \;\; \bdf_0(\bdx_1) \\
\begin{bmatrix} J_0(\bdx_1) \\ R_1 \end{bmatrix} \bdx_2 -
\left[ \begin{array}{l} \bdo \\ \bde_1 \end{array} \right]
\end{array} \right]$}
 ~~\equiv~~ \mbox{\scriptsize $
\begin{bmatrix}  \bdf_0(\bdx_1) \\ \dd{\bdx_2} \bdf_0(\bdx_1) \\ R_1 \bdx_2
-\bde_1
\end{bmatrix}$},
\end{equation}
where ~$R_1$ ~is a random matrix whose row dimension equals to the
nullity of ~$J_0(\bdx_1)$.
~The values of ~$\bdx_2 = \hbdx_2 \neq \bdo$ ~induce a
functional ~$\dd{\hbdx_2}[\hbdx_1] \in
\cD_{\hbdx}(\bdf)$.
~If the zero ~$(\hbdx_1, \hbdx_2)$ ~of ~$\bdf_1$ ~remains multiple,
then the Jacobian
~$J_1(\hbdx_1,\hbdx_2)$~ of ~$\bdf_1(\bdx_1, \bdx_2)$ ~at
~$(\hbdx_1,\hbdx_2)$ ~has a nullity ~$k_1 > 0$ ~and a nontrivial kernel.
~The depth-deflation process can be applied to ~$\bdf_1$ ~the same way as
(\ref{lvzdfsys2}) applied to $\bdf_0$.
~Namely, we seek a zero ~$(\hbdx_1,\hbdx_2,\hbdx_3,\hbdx_4)$~
to the system
\[ \bdf_2(\bdx_1,\bdx_2,\bdx_3,\bdx_4) ~~=~~
\mbox{\scriptsize $
\left[ \begin{array}{c} \;\; \bdf_1(\bdx_1, \bdx_2) \\
\begin{bmatrix} J_1( \bdx_1, \bdx_2 ) \\ R_2 \end{bmatrix}
\begin{bmatrix} \bdx_3 \\ \bdx_4 \end{bmatrix} -
\begin{bmatrix} \bdo \\ \bde_1 \end{bmatrix}
\end{array} \right]$}
\]
where ~$R_2$ ~is any matrix of size ~$k_1\times 2s$~ that makes
~$\begin{bmatrix} J_1( \bdx_1, \bdx_2 ) \\ R_2 \end{bmatrix}$~ full rank.
~By (\ref{nabla}) -- (\ref{drep}),
~equation ~$J_1 \mbox{\scriptsize $\begin{pmatrix}
\bdx_1, \bdx_2 \end{pmatrix}\begin{bmatrix} \bdx_3 \\ \bdx_4
\end{bmatrix}$} = \bdo$~ implies
\begin{equation} \label{2nddf} \begin{array}{l}
\mbox{\scriptsize $
\left[ \begin{array}{lcl}
(\bdx_3 \cdot \Dl_{\bdx_1}) \bdf_0(\bdx_1) & + &
(\bdx_4 \cdot \Dl_{\bdx_2} )\bdf_0(\bdx_1) \\
(\bdx_3 \cdot \Dl_{\bdx_1}) \dd{\bdx_2} \bdf_0(\bdx_1) & + &
(\bdx_4 \cdot \Dl_{\bdx_2}) \dd{\bdx_2} \bdf_0(\bdx_1) \\
(\bdx_3 \cdot \Dl_{\bdx_1}) (R_1 \bdx_2-\bde_1) & + &
(\bdx_4 \cdot \Dl_{\bdx_2}) (R_1 \bdx_2-\bde_1)
\end{array} \right]$}  
~~=~~
\mbox{\scriptsize $
\left[ \begin{array}{r}
\dd{\bdx_3} \bdf_0(\bdx_1) \\ (\dd{\bdx_3} \dd{\bdx_2} +\dd{\bdx_4})
\bdf_0(\bdx_1)
\\ R_1 \bdx_4 \;\;\;\;\;\;\;\;\;\;\; \end{array} \right]$} 
\;=\; \bdo  .
\end{array}
\end{equation}
Thus, the second depth-deflation seeks a solution
~$(\hbdx_1,\hbdx_2,\hbdx_3,\hbdx_4)$ ~to equations
\begin{equation} \label{df2}
 \bdf_0(\bdx_1) = \bdo, \;\;\;\; \dd{\bdx_2} \bdf_0(\bdx_1) = \bdo, \;\;\;\;
\dd{\bdx_3} \bdf_0(\bdx_1) = \bdo,  \;\;\;\;
(\dd{\bdx_3} \dd{\bdx_2} + \dd{\bdx_4}) \bdf_0(\bdx_1) = \bdo.
\end{equation}
It is important to note that ~$\hbdx_3 \ne \bdo$.
~Otherwise, from (\ref{2nddf})
\[
\mbox{\scriptsize $
\begin{bmatrix} \dd{\hbdx_4} \bdf_0(\hbdx_1) \\ R_1 \hbdx_4 \end{bmatrix}$}
~~\equiv~~
\mbox{\scriptsize $
\begin{bmatrix} J_0(\hbdx_1) \\ R_1 \end{bmatrix}$} \hbdx_4 ~~=~~ \bdo ,
\]
which would lead to ~$\hbdx_4 = \bdo$, ~making it impossible
for ~$R_2 \mbox{\footnotesize $\begin{bmatrix} \hbdx_3 \\ \hbdx_4
\end{bmatrix}$} = \bde_1$.

After ~$\al$~ depth-deflation steps, in general, we have an isolated zero
~$(\hbdx_1,\cdots,\hbdx_{2^\al})$~ to the expanded system
~$\bdf_\al (\bdx_1,\cdots,\bdx_{2^\al})$ ~with Jacobian
~$J_\al (\bdx_1,\cdots,\bdx_{2^\al})$ ~of rank ~$r_\al$.
~If ~$r_\al < 2^\al s$, ~then the next depth-deflation step seeks a zero to
~$\bdf_{\al+1}(\bdx_1,\cdots,\bdx_{2^{\al+1}}) = \bdo$
~defined in (\ref{dflta}).

\begin{lem}  \label{l:nez}
~Let ~$\bdf_0(\bdx_1) \equiv \bdf(\bdx)$~ be a system of ~$t$ ~functions of
~$s$ ~variables with a multiple zero ~$\hbdx_1 = \hbdx$.
~Assume the depth-deflation process described above reaches the extended
system ~$\bdf_{\al+1}$~ in {\em (\ref{dflta})} with isolated zero
~$(\hbdx_1,\cdots,\hbdx_{2^{\al+1}})$.
~Then ~$\hbdx_{2^j+1} \neq \bdo, \;\;\; j = 0, 1, \cdots, \al$. 
\end{lem}

\prf
~The assertion is true for ~$j=0$ ~and ~$j=1$ ~as shown above.
~Let
\[
\bdy = \mbox{\scriptsize $\left[ \begin{array}{l} \bdx_{1} \\ \;\vdots \\
\bdx_{2^{\al-1}} \end{array} \right]$}, \;\;
\bdz = \mbox{\scriptsize $
\left[ \begin{array}{l} \bdx_{2^{\al-1}+1} \\ \;\;\vdots \\
\bdx_{2^{\al-1}+2^{\al-1}} \end{array} \right]$}, \;\;
\bdu = \mbox{\scriptsize $
\left[ \begin{array}{l} \bdx_{2^\al+1} \\ \;\;\vdots \\
\bdx_{2^\al+2^{\al-1}} \end{array} \right]$}, \;\;
\bdv = \mbox{\scriptsize $
\left[ \begin{array}{l} \bdx_{2^\al+2^{\al-1}+1} \\ \;\;\vdots \\
\bdx_{2^\al+2^{\al-1}+2^{\al-1}} \end{array} \right]$}.
\]
Then
\begin{equation} \label{jal}
 J_\al (\bdy,\bdz) \mbox{\scriptsize $
\begin{bmatrix} \bdu \\ \bdv \end{bmatrix}$}
~~=~~  \mbox{\scriptsize $
\left[ \begin{array}{r} \bdu \cdot \Dl_\bdy \bdf_{\al\mns 1}(\bdy) \\
\left[
(\bdu \cdot \Dl_\bdy)(\bdz \cdot \Dl_\bdy) + (\bdv \cdot \Dl_\bdy)
\right]
\bdf_{\al\mns 1}(\bdy) \\
R_{\al\mns 1} \bdv \end{array}\right]$} ~~=~~ \bdo
\end{equation}
together with ~$\bdu = \bdo$ ~would imply
\[
J_\al (\hbdy,\hbdz) \mbox{\scriptsize $
\begin{bmatrix} \bdo \\ \bdv \end{bmatrix}$}
~~=~~  \mbox{\scriptsize $\left[ \begin{array}{r} \bdo \\
(\bdv \cdot \Dl_\hbdy) \bdf_{\al\mns 1}(\hbdy) \\
R_{\al\mns 1} \bdv \end{array}\right]$} ~~=~~
\mbox{\scriptsize $\left[ \begin{array}{l} \bdo \\
J_{\al\mns 1}(\hbdy) \\
R_{\al\mns 1} \end{array}\right]$} \bdv ~~=~~ \bdo
\]
and thereby ~$\bdv = \bdo$ ~since
~{\footnotesize $\begin{bmatrix} J_{\al\mns 1}(\hbdy) \\
R_{\al\mns 1} \end{bmatrix}$} ~is of full column rank.
~Therefore
\begin{equation} \label{hbdunez}
\hbdu ~~=~~ \left( \hbdx_{2^\al+1}^\top,\cdots,
\hbdx_{2^\al+2^{\al-1}}^\top \right)^\top ~~\ne~~ \bdo.
\end{equation}
Moreover, from (\ref{jal})
\begin{equation} \label{hbduDl}
\bdo ~~=~~ \hbdu \cdot \Dl_\bdy \bdf_{\al\mns 1}(\hbdy) ~~\equiv~~
J_{\al\mns 1}(\hbdy) \hbdu.
\end{equation}
It now suffices to show that for all ~$\eta$,
\begin{equation} \label{jeta}
 J_\eta (\hbdx_1,\cdots,\hbdx_{2^\eta})
\mbox{\scriptsize $
\begin{bmatrix} \bdw_1 \\ \vdots \\ \bdw_{2^\eta} \end{bmatrix}$}
~~=~~ \bdo ~~\mbox{ \ \ and \ \ }~~
\mbox{\scriptsize $
\begin{bmatrix} \bdw_1 \\ \vdots \\ \bdw_{2^\eta} \end{bmatrix}$}
~~\ne~~ \bdo
\end{equation}
would imply ~$\bdw_1 \ne \bdo$.
~Obviously, this is true for ~$\eta = 1$.
~Assume it is true up to ~$\eta-1$.
~Then, using the same argument for (\ref{hbdunez}) and  (\ref{hbduDl}),
we have (\ref{jeta}) implying
\[
\mbox{\scriptsize $
\left[ \begin{array}{l} \bdw_1 \\ \;\vdots \\ \bdw_{2^{\eta-1}}
\end{array} \right]$}
~~\ne~~ \bdo \mbox{ \ \ and \ \ }
J_{\eta-1} \mbox{\scriptsize $
\left[ \begin{array}{l} \bdw_1 \\ \;\vdots \\ \bdw_{2^{\eta-1}}
\end{array} \right]$}
~~=~~ \bdo.
\]
Thus ~$\bdw_1 \ne \bdo$~ from the induction assumption. \hfill \qed

It is clear that the third depth-deflation, if necessary, adds variables
~$\bdx_5$, ~$\bdx_6$, ~$\bdx_7$, ~$\bdx_8$~ and equations
\begin{equation}
\begin{array}{l}
\dd{\bdx_5} \bdf(\bdx_1) = \bdo, \;\;\;\;
(\dd{\bdx_5} \dd{\bdx_2} + \dd{\bdx_6}) \bdf(\bdx_1) = \bdo,
\;\;\;\;\;\;\;
(\dd{\bdx_5} \dd{\bdx_3} + \dd{\bdx_7}) \bdf(\bdx_1) = \bdo, \\
(\dd{\bdx_5} \dd{\bdx_3} \dd{\bdx_2} +
\dd{\bdx_5} \dd{\bdx_4} +
\dd{\bdx_3} \dd{\bdx_6} + \dd{\bdx_7} \dd{\bdx_2} + \dd{\bdx_8})
\bdf(\bdx_1) = \bdo.
\end{array} \label{df3}
\end{equation}
Any solution ~$(\hbdx_1,\cdots,\hbdx_8) \in \C^{8s}$~ to (\ref{df2})
and (\ref{df3}) induces eight differential functionals
\[ \begin{array}{l}
1, \;\;\; \dd{\hbdx_2}, \;\;\;\;\; \dd{\hbdx_3}, \;\;\;\;\; \dd{\hbdx_5},
~~\dd{\hbdx_3} \dd{\hbdx_2} + \dd{\hbdx_4}, \;\;\;\;
\dd{\hbdx_5} \dd{\hbdx_2} + \dd{\hbdx_6}, \;\;\;\;
\dd{\hbdx_5} \dd{\hbdx_3} + \dd{\hbdx_7}, \\
\dd{\hbdx_5} \dd{\hbdx_3} \dd{\hbdx_2} + \dd{\hbdx_5} \dd{\hbdx_4}
+ \dd{\hbdx_3} \dd{\hbdx_6} + \dd{\hbdx_7} \dd{\hbdx_2} + \dd{\hbdx_8}
\end{array} \]
that vanish on ~$\bdf$ ~at ~$\hbdx_1$.
~In general, the ~$\al$-th depth-deflation step produces a collection of
~$2^\al$ ~differential functionals of order ~$\al$ ~or less that vanish on the
system ~$\bdf$ ~at ~$\hbdx_1$.
~Also notice that the highest order differential terms are
\[ \dd{\hbdx_2} \equiv \dd{\hbdx_{2^0+1}}, \;\;
\dd{\hbdx_3}\dd{\hbdx_2} \equiv \dd{\hbdx_{2^1+1}} \dd{\hbdx_{2^0+1}},
\;\; \dd{\hbdx_5}\dd{\hbdx_3}\dd{\hbdx_2} \equiv
\dd{\hbdx_{2^2+1}}\dd{\hbdx_{2^1+1}}\dd{\hbdx_{2^0+1}} \]
for depth-deflation steps 1, 2 and 3, respectively.

Actually these functionals induced by the depth-deflation method all
belong to the dual space ~$\cD_\hbdx(\bdf)$.
~To show this, we define differential operators
~$\Phi_\al$, $\al = 1,2,\cdots$ ~as follows.
\begin{equation}
\Phi_{\nu+1} \; =  \; \sum_{\zeta = 1}^{2^\nu}
\; \bdx_{_{2^\nu+\zeta}}\cdot \Dl_{\bdx_\zeta}, \;\;\;
\nu = 0, 1, \cdots. \label{Phi}
\end{equation}
Specifically, ~$\Phi_1 \,=\, \bdx_2 \cdot \Dl_{\bdx_1}$,  
~$\Phi_2 \,=\, \bdx_3 \cdot \Dl_{\bdx_1} + \bdx_4 \cdot \Dl_{\bdx_2}$
~and 
~$\Phi_3 \,=\, \bdx_5 \cdot \Dl_{\bdx_1} + \bdx_6 \cdot \Dl_{\bdx_2}
+ \bdx_7 \cdot \Dl_{\bdx_3} + \bdx_8 \cdot \Dl_{\bdx_4}$.
~For convenience, let ~$\Phi_0$~ represent the identity operator.
~Thus
\[ \begin{array}{l}
\Phi_0 \bdf(\bdx_1) ~=~ \bdf(\bdx_1), ~~~~
\Phi_1 \bdf(\bdx_1) ~=~ \dd{\bdx_2} \bdf (\bdx_1), ~~~~
\Phi_2 \bdf(\bdx_1) ~=~ \dd{\bdx_3} \bdf (\bdx_1), \\
\Phi_2 \circ \Phi_1 \bdf(\bdx_1)  ~~=~~
(\bdx_3 \cdot \Dl_{\bdx_1}) \dd{\bdx_2} \bdf (\bdx_1) +
(\bdx_4 \cdot \Dl_{\bdx_2}) \dd{\bdx_2} \bdf (\bdx_1) ~~=~~
(\dd{\bdx_3} \dd{\bdx_2} + \dd{\bdx_4}) \bdf(\bdx_1)
\end{array} \]
etc.
~For any expanded system ~$\bdf_\al(\bdx_1,\cdots,\bdx_{2^\al})$~ generated
in the depth-deflation process, its Jacobian
~$J_\al(\bdx_1,\cdots,\bdx_{2^\al})$~satisfies
\[  J_\al(\bdx_1,\cdots,\bdx_{2^\al})
\mbox{\scriptsize $
\left[ \begin{array}{l} \bdx_{2^\al + 1} \\ \;\;\;\vdots \\
\bdx_{2^\al + 2^\al} \end{array} \right]$}
~~=~~ \Phi_{\al+1} \bdf_\al(\bdx_1, \cdots,\bdx_{2^\al}).
\]
~It is easy to see that
(\ref{df2}) and (\ref{df3}) can be written as
\[ \begin{array}{l}
\Phi_0 \bdf(\bdx_1)  ~=~ \bdo, ~~~~
\Phi_1 \bdf(\bdx_1)  ~=~ \bdo, ~~~~
\Phi_2 \bdf(\bdx_1)  ~=~  \bdo, \;\;\;\; \Phi_2 \circ \Phi_1 \bdf(\bdx_1)
~~=~~  \bdo, \\
\Phi_3 \bdf(\bdx_1) ~~=~~\bdo, \;\;\;\; \Phi_3 \circ \Phi_1 \bdf(\bdx_1)
~~=~~ \bdo,
\;\;\;\;\;
\Phi_3 \circ \Phi_2 \bdf(\bdx_1)  ~~=~~\bdo, \; \;\;\;\;
\Phi_3 \circ \Phi_2 \circ \Phi_1 \bdf(\bdx_1) ~~=~~ \bdo.  \end{array}
\]

As a consequence, Theorem \ref{t:conj1} given in \S\ref{s:dm} provides
an upper bound, {\em the depth}, on the
number of depth-deflation steps required to regularize the singularity
at the multiple zero.
~This bound substantially
improves the result in \cite[Theorem 3.1]{lvz06}.
~In fact, our version of the deflation method deflates {\em depth}
rather than the multiplicity as suggested in \cite{lvz06}.

%
{\bf Proof of Theorem \ref{t:conj1}}.
~We first claim that the ~$\al$-th depth-deflation step 
induces {\em all} differential functionals 
\begin{equation} \label{phifl}
f ~~\longrightarrow~~ 
\Phi_{\mu_1}\circ \cdots \circ \Phi_{\mu_k}\,f
\big|_{(\bdx_1,\cdots,\bdx_{2^\al}) = (\hbdx_1,\cdots,\hbdx_{2^\al})}
~~~\mbox{with}~~
\al \geq \mu_1 > \mu_2 > \cdots > \mu_k \geq 0
\end{equation}
and ~$1 \leq k \leq \al$ ~that vanish on ~$\bdf$.
~This is clearly true for ~$\al=1$ ~since ~$\bdf_1(\bdx_1,\bdx_2)=\bdo$
~induces ~$\Phi_0 \bdf(\bdx_1) \,=\, \Phi_1 \bdf(\bdx_1) \,\equiv\, 
\Phi_1 \Phi_0 f(\bdx_1) \,=\,\bdo$ ~at ~$(\bdx_1,\bdx_2) = (\hbdx_1,\hbdx_2)$.
~Assume the claim is true for ~$\al-1$.
~At the ~$\al$-th depth-deflation, consider a functional (\ref{phifl}).
~If ~$\mu_1 < \al$, ~then such a functional has already been
induced from solving
~$\bdf_{\al-1} = \bdo$.
~On the other hand, if ~$\mu_1 = \al$, ~then
~$\Phi_{\mu_2}\circ \cdots \circ \Phi_{\mu_k} \bdf(\bdx_1) \,=\, \bdo$,
~for ~$\al-1 \geq \mu_2 > \cdots > \mu_k \geq 0$
~is in ~$\bdf_{\al-1} = \bdo$.
~Therefore ~$\Phi_\al \bdf_{\al-1}$~ induces the functional in (\ref{phifl}).
~Next, the functional in (\ref{phifl}) satisfies
closedness condition (\ref{cxjfi}).
~To show this, let ~$p$ ~be any polynomial in variables ~$\bdx$.
~By applying the product rule
~$ \Phi_\al (f\,g) = (\Phi_\al\, f)\,g + (\Phi_\al\, g)\, f$~
in an induction,
\[ \Phi_{\mu_1}\circ \cdots \circ \Phi_{\mu_k} (p f_i)
= \sum_{\{\eta_1, \cdots, \eta_j\} \subset \{\mu_1,\cdots,\mu_k\} }
p_{\eta_1 \cdots \eta_j} \Phi_{\eta_1} \circ \cdots \circ \Phi_{\eta_j} f_i
\]
where ~$\eta_1 > \cdots > \eta_j$~ and ~$p_{\eta_1 \cdots \eta_j}$ ~is a
polynomial generated by applying ~$\Phi_j$'s on ~$p$.
~Therefore ~$\Phi_{\mu_1}\circ \cdots \circ \Phi_{\mu_k} (p f_i) = 0$ ~at
~$(\hbdx_1,\cdots,\hbdx_{2^\al})$ ~since
~$\Phi_{\eta_1} \circ \cdots \circ \Phi_{\eta_j} f_i = 0$,
~showing that functionals (\ref{phifl}) all belong to
~$\cD_{\hbdx}(\bdf)$.
~Finally, the highest order part of the differential functional
~$\Phi_\al \circ \Phi_{\al-1} \circ \cdots \circ \Phi_1$
~is
~$\mbox{$\prod_{j=0}^{\al-1} (\hbdx_{2^j+1}\cdot \Dl_\bdx) ~~\equiv~~
\prod_{j=0}^{\al-1} \dd{\hbdx_{2^j+1}}$}$ 
~which is of order ~$\al$ ~since ~$\hbdx_{2^j+1} \ne \bdo$ ~by
Lemma~\ref{l:nez}.
~However, differential orders of all functionals in ~$\cD_{\hbdx}(\bdf)$
~are bounded by ~$\dl_{\hbdx}(\bdf)$, ~so is ~$\al$.
\hfill \qed

In general, Theorem \ref{t:conj1} does not guarantee those
~$2^k$ ~functionals are linearly independent.
~From computing experiments, the number ~$k$ ~of depth-deflation
steps also correlates to the breadth
$\bt_{\hbdx}(\bdf)$.
~Especially when ~$\bt_{\hbdx}(\bdf) = 1$, ~it appears that
~$k$ ~always reaches its maximum.
~This motivates the special case breadth-one algorithm
which will be presented in \S\ref{s:b1}.
~On the other hand, when breadth ~$\bt_{\hbdx}(\bdf) > 1$,
~very frequently the depth-deflation process pleasantly terminates
only after {\em one} depth-deflation step regardless of the
depth or multiplicity.
~A possible explanation for such a phenomenon is as follows.
~At each depth-deflation step, say the first, the isolated zero 
~$\hbdz$ ~to the system (\ref{lvzdfsys2}) is multiple only if there
is a differential functional in the form of
~$\dd{\bdx_3} \dd{\bdx_2} + \dd{\bdx_4}$~ in ~$\cD^2_\hbdx (\bdf)$
~while ~$R_1 \bdx_2 = \bde_1$~ and ~$R_1 \bdx_4 = \bdo$~ for a randomly
chosen ~$R_1$.
~In most of the polynomial systems we have tested, functionals in
this special form {\em rarely} exist in ~$\cD^2_\hbdx (\bdf)$
~when ~$\bt_{\hbdx}(\bdf) > 1$.
~If no such functionals exist in ~$\cD^2_{\hbdx}(\bdf)$, ~the zero
~$\hbdz$ ~must be a simple zero of ~$\tilde{F}$ ~in (\ref{lvzdfsys2})
according to Theorem \ref{t:conj1},
therefore the depth-deflation ends at ~$k=1$~ step.

\vspace{-4mm}
\subsection{Special case: dual space of breadth one} \label{s:b1}
\vspace{-4mm}

Consider a nonlinear system
~$\bdf = [ f_1,\cdots,f_t ]^\top$~ having breadth one at an
isolated zero ~$\hbdx$, ~namely ~$\bt_{\hbdx}(\bdf) = 1$.
~The Hilbert function is
~$\{1,1,\cdots,1,0,\cdots \}$,
~making the depth ~one less than the multiplicity:
~$\dl_{\hbdx}(\bdf) = \dm\big(\cD_\hbdx(\bdf)\big)-1$.
~This special case includes the most fundamental
univariate equation ~$f(x) ~=~ 0$ ~at a multiple zero.
~As mentioned above, the general depth-deflation method
derived in \S\ref{s:dm} {\em always}
exhausts the maximal number of steps in this case,
and the final system is expanded undesirably from ~$t \times s$
~to over ~$(2^{m-1} t) \times (2^{m-1} s)$ ~at an ~$m$-fold zero.
~To overcome this exponential growth of the system size, we
shall modify the depth-deflation process for breadth-one system
in this section so that
the regularized system is of size close to ~$(m t) \times (m s)$,
~and upon solving the system, a complete basis
for the dual space ~$\cD_\hbdx(\bdf)$ ~is obtained as a by-product.

Denote ~$\bdx = \bdx_1$ ~and the zero ~$\hbdx = \hbdx_1$ ~as in \S\ref{s:dm}.
~It follows from (\ref{S1}), the breadth
~$\bt_{\hbdx}(\bdf) \; = \; \rmh(1) =  
\nullity{J_0(\hbdx_1)}  =  1 $~
implies system (\ref{lvzdfsys2}), simplifying to
~$\mbox{\footnotesize $\begin{bmatrix}  J_0(\hbdx_1) \\
\bdb^\h \end{bmatrix}$}
\bdx_2 = \mbox{\footnotesize $\begin{bmatrix} \bdo \\ 1
\end{bmatrix}$}$~ in the variable vector ~$\bdx_2$,
~has a unique solution ~$\hbdx_2 \in \C^s$~
for randomly chosen vector ~$\bdb \in \C^s$.
~Similar to the general depth-deflation method in \S~\ref{s:dm}, the first step
of depth-deflation is to expanded the system:
\begin{eqnarray} \label{G1F1}
&& \bdg_1 \left( \bdx_1 , \bdx_2  \right)
~~=~~ \mbox{\scriptsize $
\left[ \begin{array}{l} \bdh_0(\bdx_1) \\ \bdh_1(\bdx_1,\bdx_2)
\end{array} \right]$} \\
&& \mbox{ \ where \ }
\bdh_0 (\bdx_1) ~\equiv~ \bdf(\bdx) ~~\mbox{and}~~
\bdh_1(\bdx_1 ,\bdx_2) ~=~
\mbox{\scriptsize $
\begin{bmatrix}J_0(\bdx_1) \,\bdx_2  \\ \bdb^\h \bdx_2 -1
\end{bmatrix}$} ~\equiv~
\mbox{\scriptsize $
\begin{bmatrix}  \dd{\bdx_2} \bdf(\bdx_1)  \\ \bdb^\h \bdx_2 -1
\end{bmatrix}$}. \nonumber
\end{eqnarray}
The system ~$\bdg_1(\bdx_1,\bdx_2)$ ~has an isolated zero
~$(\hbdx_1, \hbdx_2)$.
~If the Jacobian ~$J_1(\bdx_1,\bdx_2)$~ of ~$g_1(\bdx_1,\bdx_2)$~ is of full
rank at $(\hbdx_1,\hbdx_2)$, ~then the system is regularized and the 
depth-deflation process terminates.
~Otherwise, there is a nonzero vector $(\bdv_1,\bdv_2) \in \C^{2s}$
~such that
\begin{equation}  \label{J1}
J_1(\hbdx_1,\hbdx_2) \mbox{\scriptsize $
\left[ \begin{array}{l} \bdv_1 \\ \bdv_2 \end{array}
\right]$} ~~\equiv~~
\mbox{\scriptsize $
\left[ \begin{array}{l} \dd{\bdv_1} \bdf(\hbdx_1) \\
(\dd{\bdv_1} \dd{\hbdx_2} + \dd{\bdv_2} ) \bdf(\hbdx_1)\\
\bdb^\h \bdv_2
\end{array} \right]$} ~~=~~ \bdo.
\end{equation}
Since the Jacobian ~$J_0(\hbdx)$ ~of ~$\bdf$ ~at ~$\hbdx_1$ is of nullity one,
there is a constant ~$\gamma \in \C$ ~such that ~$\bdv_1  = \gamma \hbdx_2$.
~Equation (\ref{J1}) together with
~$\bt_{\hbdx_0}(\bdf) = 1$~ and ~$(\bdv_1, \bdv_2) \neq (\bdo, \bdo)$~
imply ~$\gamma \neq 0$.
~Consequently we may choose ~$\gamma = 1$, ~namely ~$\bdv_1 = \hbdx_2$.
~Setting ~$\hbdx_3 = \bdv_2$, ~the system
\begin{eqnarray} \label{G2F2}
\bdg_2(\bdx_1,\bdx_2,\bdx_3) & \equiv  &
\mbox{\scriptsize $
\left[ \begin{array}{l} \bdh_0(\bdx_1) \\ \bdh_1(\bdx_1,\bdx_2) \\
\bdh_2(\bdx_1,\bdx_2,\bdx_3) \end{array} \right]$}
~~=~~ \mbox{\scriptsize $\left[ \begin{array}{r} \bdf(\bdx_1) \\
\dd{\bdx_2} \bdf(\bdx_1)  \\ \bdb^\h \bdx_2 -1 \\
(\dd{\bdx_2} \dd{\bdx_2} + \dd{\bdx_3})\bdf(\bdx_1)
\\ \bdb^\h \bdx_3 \end{array} \right]$}
\\ & &
\mbox{ \ \ where \ }
\bdh_2(\bdx_1,\bdx_2,\bdx_3) ~~=~~ \mbox{\scriptsize $\begin{bmatrix}
(\dd{\bdx_2} \dd{\bdx_2} + \dd{\bdx_3})\bdf(\bdx_1)
\\ \bdb^\h \bdx_3
\end{bmatrix}$}  \nonumber
\end{eqnarray}
has an isolated zero ~$(\hbdx_1, \hbdx_2, \hbdx_3)$.
~In general,
if an isolated zero ~$(\hbdx_1,\cdots,\hbdx_{\gamma+1})$ ~to
the system
\[ \bdg_{\gamma} (\bdx_1,\cdots,\bdx_{\gamma+1})
~~=~~ \mbox{\scriptsize $
\left[ \begin{array}{l} \bdh_0(\bdx_1) \\ \bdh_1(\bdx_1,\bdx_2)
\\ ~~~\vdots \\ \bdh_{\gamma}(\bdx_1,\cdots,\bdx_{\gamma+1}) \end{array} 
\right]$}
\]
remains singular, or the Jacobian
~$J_{\gamma}(\hbdx_1,\cdots,\hbdx_{\gamma+1})$~ is rank-deficient, then
there is a non-zero solution to the homogeneous system
\[ J_{\gamma}(\hbdx_1,\cdots,\hbdx_{\gamma+1}) 
\mbox{\scriptsize $\left[ \begin{array}{l}
\bdu_1 \\ ~~\vdots \\ \bdu_{\gamma +1}\end{array} \right]$}
~~\equiv~~ \mbox{\scriptsize $\left[ \begin{array}{c}
J_{\gamma-1}(\hbdx_1,\cdots,\hbdx_{\gamma}) \left[ \begin{array}{l}
\bdu_1 \\ ~\vdots \\ \bdu_{\gamma}\end{array} \right] \\ * \end{array}
\right]$} ~~=~~ \bdo.
\]
Therefore, by setting ~$\bdu_j = \hbdx_{j\pls 1}$ ~for ~$j=1,\ldots,\gamma$,
~we take its unique solution
~$\bdu_{\gamma+1}$ ~as ~$\hbdx_{\gamma+2}$.

The pattern of this depth-deflation process can be illustrated by defining
\begin{equation} \label{Psi}
\Psi  ~~=~~  \sum_{\eta = 1}^\infty \bdx_{\eta+1} \cdot \Delta_{\bdx_{\eta}}.
\end{equation}
When applying ~$\Psi$ ~to any function ~$f$ ~in (vector)
variables, say  ~$\bdx_1, \cdots, \bdx_\sg$, the resulting
~$\Psi f$~ is a {\em finite} sum
since ~$\Dl_{\bdx_\mu} f = \bdo$~ for ~$\mu \ge \sg+1$.
~Thus,
\begin{eqnarray}
\lefteqn{\bdh_1 (\bdx_1, \bdx_2) ~~=~~
\mbox{\scriptsize $\begin{bmatrix}
\Psi \bdh_0(\bdx_1) \\  \bdb^\h \bdx_2 -1 \end{bmatrix}$},
~~~\bdh_2(\bdx_1,\bdx_2,\bdx_3) ~~=~~
\mbox{\scriptsize $\begin{bmatrix}
\Psi \bdh_1(\bdx_1,\bdx_2) \\  \bdb^\h \bdx_3 -1 \end{bmatrix}$}
~~~~\mbox{and}} \nonumber \\
& & \nonumber \\
\label{Fnu}
& & \bdh_\nu(\bdx_1,\cdots,\bdx_{\nu})  ~~=~~
\mbox{\scriptsize $\left[ \begin{array}{c}
\mbox{$ \overbrace{\Psi \circ \Psi \circ \cdots \circ \Psi}^{\nu-1} $}
\,\bdh_1 (\bdx_1,\bdx_2), \\
\\ \bdb^\h \bdx_{\nu\pls 1} \end{array} \right]$},
~~~~\mbox{for}~~~ \nu \geq 2.
\end{eqnarray}
For instance, with ~$\bdh_1$ ~and ~$\bdh_2$ ~in (\ref{G1F1}) and
(\ref{G2F2}) respectively, we have
\[ \bdh_3(\bdx_1, \bdx_2, \bdx_3, \bdx_4 ) ~~=~~ \mbox{\scriptsize
$\begin{bmatrix}
(\dd{\bdx_2} \dd{\bdx_2} \dd{\bdx_2} + 3 \dd{\bdx_2} \dd{\bdx_3}
+ \dd{\bdx_4}) \bdh_0(\bdx_1) \\
\bdb^\h \bdx_4
\end{bmatrix}$}.
\]
If, say, ~$\bdh_3 = \bdo$~ at ~$(\hbdx_1, \hbdx_2, \hbdx_3, \hbdx_4 )$,
~a functional ~$f \longrightarrow
 \left(\dd{\hbdx_2} \dd{\hbdx_2} \dd{\hbdx_2} + 3 \dd{\hbdx_2} \dd{\hbdx_3}
+ \dd{\hbdx_4}\right)f(\bdx_1)$
~is obtained and it vanishes on the system ~$\bdf$.
~The original system ~$\bdf(\bdx)=\bdo$ ~provides a trivial functional
~$\rd_{0\cdots 0} \;:\; f \rightarrow f(\hbdx_1)$.
~By the following lemma those functionals are all in the dual space.

\begin{lem} \label{l:b1func}
~Let ~$\bdf = [f_1,\cdots,f_t]^\top$ ~be a nonlinear system
with an isolated zero ~$\hbdx \in \C^s$.
~Write ~$\bdg_0 = \bdf$, ~$\hbdx_1 = \hbdx$ ~and ~$\bdx_1 = \bdx$.
~For any ~$\gamma \in \{\,1,2,\cdots\,\}$, ~let ~$(\hbdx_1,\, \hbdx_2,\,
\cdots, \, \hbdx_{\gamma\pls 1})$~ be a zero of
\begin{equation} \label{Ggamma}
 \bdg_\gamma (\bdx_1,\bdx_2,\cdots, \bdx_{\gamma\pls 1} )
\; =\; \mbox{\scriptsize $ \left[ \begin{array}{l}
\bdh_0(\bdx_1) \\ \;\; \vdots \;\;\;\;\;\; \ddots \\
\bdh_\gamma (\bdx_1,\cdots,\bdx_{\gamma+1}) \end{array} \right]$}.
\end{equation}
Then the functionals derived from
~$\bdg_\gamma (\hbdx_1,\cdots,\hbdx_{\gamma+1}) = \bdo$
~constitutes a linearly independent subset of the dual space
~$\cD_{\hbdx_0}(\bdf)$.
\end{lem}

\prf
By a rearrangement, finding a zero of 
~$\bdg_\gamma (\bdx_1,\bdx_2,\cdots, \bdx_{\gamma\pls 1} )$~ is equivalent
to solving
\begin{equation} \label{fbdx1}
\begin{array}{rcl}
\bdf(\bdx_1) = \bdo, && \bdb^\h \bdx_{2} = 1, \\
\Psi \bdf (\bdx_1) = \bdo, && \bdb^\h \bdx_3 = 0, \\
\vdots \;\;\;\;\;  && ~~~\vdots \\
\Psi \circ \cdots \circ \Psi \bdf(\bdx_1) = \bdo, && \bdb^\h
\bdx_{\gamma\pls 1} = 0.
\end{array}
\end{equation}
for ~$(\bdx_1,\cdots,\bdx_{\gamma\pls 1}) \in \C^{(\gamma\pls 1)s}$.
~Let ~$(\hbdx_1,\cdots, \hbdx_{\gamma+1})$~ be an isolated zero.
~Then each ~$\Psi \circ \cdots \circ \Psi$~ induces a differential
functional
\begin{equation}
\rho_\al ~~:~~
f \longrightarrow \overbrace{\Psi \circ \cdots \circ \Psi}^\al f
\bigg|_{
(\bdx_1,\cdots,\bdx_{\al+1}) =
(\hbdx_1,\cdots, \hbdx_{\al+1})},  \label{rhoal}
~~~\mbox{for}~~~ \al = 0, 1, \cdots, \gamma.
\end{equation}
Those functionals vanish on ~$f_1,\cdots,f_t$~ because of (\ref{fbdx1}).
~Since ~$\Psi$~ satisfies product rule
~$\Psi (fg) = (\Psi f)g+f(\Psi g)$~ for any functions ~$f$ ~and ~$g$ ~in
finitely many variables among ~$\bdx_1,\bdx_2,\cdots,$
for any polynomial ~$p \in \C[\bdx_1]$, ~we have,
for ~$\al = 0, 1, \cdots, \gamma$ ~and ~$i = 1, \cdots, t$,
\[
\rho_\al (pf_i) ~~=~~
\sum_{j=0}^\al \begin{pmatrix} \al \\ j \end{pmatrix}
(\overbrace{\Psi \circ \cdots \circ \Psi}^j p)
(\overbrace{\Psi \circ \cdots \circ \Psi}^{\al-j} f_i )
\bigg|_{ (\bdx_1,\cdots,\bdx_{\al+1}) =
(\hbdx_1,\cdots, \hbdx_{\al+1})} = 0 .
\]
Namely, ~$\rho_\al$'s satisfy the closedness condition (\ref{cxjfi}),
so they belong to ~$\cD_{\hbdx_1}(\bdf)$.

The leading (i.e., the highest order differential) term of
~$\rho_\al$~ is ~$\overbrace{\dd{\hbdx_2}\cdots \dd{\hbdx_2}}^\al$
which is of order $\al$ since ~$\hbdx_2 \ne \bdo$.
~Therefore, they are linearly independent.
\hfill \qed

\begin{thm}[Breadth-one Deflation Theorem] \label{t:b1}
~Let ~$\hbdx$~ be an isolated multiple zero of the nonlinear system
~$\bdf = [f_1,\cdots,f_t]^\top$ ~with breadth ~$\beta_{\hbdx} (\bdf) = 1$.
~Then there is an integer ~$\gamma \leq \dl_{\hbdx} (\bdf)$
~such that, for almost all ~$\bdb\in \C^s$, ~the system
~$\bdg_\gamma$ ~in {\em (\ref{Ggamma})} has a simple zero
~$(\hbdx_1,\, \hbdx_2,\, \cdots, \, \hbdx_{\gamma\pls 1})$
~which induces ~$\gamma\npls 1$~ linearly independent functionals in
~$\cD_{\hbdx}(\bdf)$.
\end{thm}

\prf A straightforward consequence of Lemma \ref{l:b1func}. \hfill \qed

While the general depth-deflation method usually terminates with one or 
two steps
of system expansion for systems of breadth higher than one, the breadth
one depth-deflation {\em always} terminates at step 
~$\gamma = \dl_\hbdx(\bdf)$ ~exactly.
~Summarizing the above elaboration, we give the pseudo-code of
an efficient algorithm for computing the
multiplicity structure of the breadth one case as follows:

\vspace{-4mm}
\begin{itemize} \parskip-0.5mm
\item[] {\bf Algorithm} {\sc BreadthOneMultiplicity}%
\item[] {\tt Input:} \ Nonlinear system ~$\bdf = [f_1,\ldots,f_t]^\h$,
~zero ~$\hbdx_1 \in \C^s$
\begin{itemize}
\item ~{\tt set} random vectors ~$\bdb \in \C^s$ ~and
{\tt obtain} ~$\hbdx_2$ ~by solving
~{\scriptsize $\left[ \begin{array}{c} J(\hbdx_1) \\ \bdb^\h
\end{array} \right]\,\bdx_2 = \left[ \begin{array}{c} \bdo \\  1
\end{array} \right]$}
\item ~{\tt initialize} ~$\bdp_2(\bdx_1,\bdx_2) = J(\bdx_1)\bdx_2$
\item ~{\tt for} ~$k=2,3, \ldots$ ~{\tt do}
\begin{itemize}
\item ~{\tt set} ~\mbox{$\bdd_k(\bdx_1,\ldots,\bdx_k) ~=~
-\sum_{j=1}^{k-1}\, \hbdx_{j\pls 1} \cdot \Dl_{\bdx_j} \,\bdp_k(\bdx_1,\ldots,
\bdx_k)$}
\item ~{\tt solve} for ~$\bdx_{k\pls 1} = \hbdx_{k\pls 1}$ ~in the system
\begin{equation} \label{Jkp1}
 \mbox{\scriptsize $\left[ \begin{array}{c} J(\hbdx_1) \\ \bdb^\h
\end{array} \right] \bdx_{k\pls 1} ~~=~~ \left[ \begin{array}{c}
\bdd_k(\hbdx_1,\ldots,\hbdx_k) \\  0 \end{array} \right]$}  \end{equation}
\item ~{\tt if} the equation (\ref{Jkp1}) has no solution, 
{\tt set} ~$\gamma = k-1$
~and \newline {\tt break} the loop; ~otherwise, {\tt set}
\[ \bdp_{k\pls 1}(\bdx_1,\ldots,\bdx_{k\pls 1}) ~=~
\Psi \,\bdp_k(\bdx_1,\ldots,\bdx_k) ~\equiv~ \bdd_k(\bdx_1,\ldots,\bdx_k)
+ J(\bdx_1)\bdx_{k\pls 1}\]
\end{itemize}
\item[] ~{\tt end do}
\end{itemize}
\item[] {\tt Output:} ~multiplicity ~$\gamma+1$ ~and functionals
~$\rho_0,~\rho_1,~\ldots,~\rho_\gamma$ ~as in (\ref{rhoal})
\end{itemize}

\begin{example}  \em ~One of the main advantages of our algorithms is the 
capability of accurate identification of multiplicity structures even if
the system data are given with perturbations and the zero is approximate.
~Consider the sequence of nonlinear systems
\end{example}
\vspace{-4mm}
\begin{equation} \label{mbc}
\tilde{\bdf}_k(x,y,z) ~~=~~ 
[\,x^2\,\sin y, ~y-z^2,~z-1.772453850905516\,\cos x^k \,]^\top,
\end{equation}
which is an inexact version of the system 
~$\bdf_k(x,y,z) = [ \,x^2\,\sin y, ~y-z^2,~z-\sqrt{\pi} \,\cos x^k ]^\top$
with breadth one and isolated zero ~$(0,\pi,\sqrt{\pi})$.
~The multiplicity is ~$2(k+1)$ ~and the depth is 
~$\dl_{(0,\pi,\sqrt{\pi})}(\bdf_k) = 2k+1$ ~for ~$k=1,2,\ldots$.
~Our code {\sc BreadthOneMultiplicity} running on floating point arithmetic
accurately identifies the multiplicity structure with the 
approximate dual basis
\begin{eqnarray*} 
& & 1, ~~\rd_{x},~~\rd_{x^2},~\ldots~,\rd_{x^{2k\mns 1}},
~~\rd_{y}+\mbox{\scriptsize 0.2820947917738781}\,\rd_z 
-\mbox{\scriptsize 0.3183098861837908}\, \rd_{x^{2k}}, \\
& & ~~~~~~~~~~\rd_{xy}+\mbox{\scriptsize 0.2820947917738781}\,\rd_{xz}
-\mbox{\scriptsize 0.3183098861837908}\, \rd_{x^{2k\pls 1}}
\end{eqnarray*}
at the numerical zero 
~$(\mbox{\scriptsize 0, 3.141592653589793, 1.772453850905516})$.
~The computing time is shown in Table~\ref{t:mbc} for 
Algorithm {\sc BreadthOneMultiplicity}.
\qed

\begin{table}[ht]
\begin{center}
\begin{tabular}{||c|r|r|r|r|r||} \hline \hline
 ~$k$: &\multicolumn{1}{|l|}{~~2}  & \multicolumn{1}{|l|}{~~4}  &
\multicolumn{1}{|l|}{~~6}  & \multicolumn{1}{|l|}{~~8}  &
\multicolumn{1}{|l||}{~~10}  \\ \hline
computed depth : & 5~~ & 9~~ & 13~~ & 17~~ & 21~~ \\ 
computed multiplicity : & 6~~ & 10~~ & 14~~ & 18~~ & 22~~ \\ 
\hline \hline
{\sc BreadthOneMultiplicity} elapsed time & ~~~~0.34 & ~~~~1.45 & ~~~~3.58 & 
~~~~18.22 & ~~~~63.42 \\
\hline \hline
\end{tabular}
\end{center}
\caption{~Results of {\sc BreadthOneMultiplicity} in floating point arithmetic
on the inexact systems ~$\tilde{\bdf}_k$ ~in (\ref{mbc}) at the
approximate zero 
$(\mbox{\scriptsize 0, 3.141592653589793, 1.772453850905516})$.}
\label{t:mbc}
\end{table}

In our extensive computing experiments,
Algorithm~{\sc BreadthOneMultiplicity} always produces a complete
dual basis without premature termination.
~We believe the following conjecture is true.

\begin{conj}
~Under the assumptions of Theorem~{\em \ref{t:b1}}, Algorithm
{\sc BreadthOneMultiplicity} terminates at ~$\gamma=\dl_\hbdx(\bdf)$
~and generates a complete basis for the dual space
\[ \cD_\hbdx(\bdf) ~~=~~ \spn \{ \rho_0, \rho_1, \ldots, \rho_\gamma \}.
\]
\end{conj}

{\bf Acknowledgements.} 
~The authors wish to thank following scholars:
~Along with many insightful discussions, 
Andrew Sommese provided a preprint \cite{bps09} which presented an 
important application of this work,
Hans Stetter provided the diploma thesis \cite{Tha96} of his former student,
~Teo Mora pointed out Macaulay's original contribution \cite{mac16} 
elaborated in his book \cite{mora2}, and
Lihong Zhi pointed out the reference \cite{lbh}.

\addcontentsline{toc}{section}{References}

\end{document}